\documentclass[12pt,reqno]{amsart}
\usepackage[a4paper,centering,total={160mm,237mm}]{geometry}

\usepackage{mathtools}
\usepackage{amssymb}
\usepackage{tikz}
\usepackage{graphicx}
\usetikzlibrary{shapes.geometric}

\newcommand\orcidicon[1]{\href{https://orcid.org/#1}{\includegraphics[scale=0.0045]{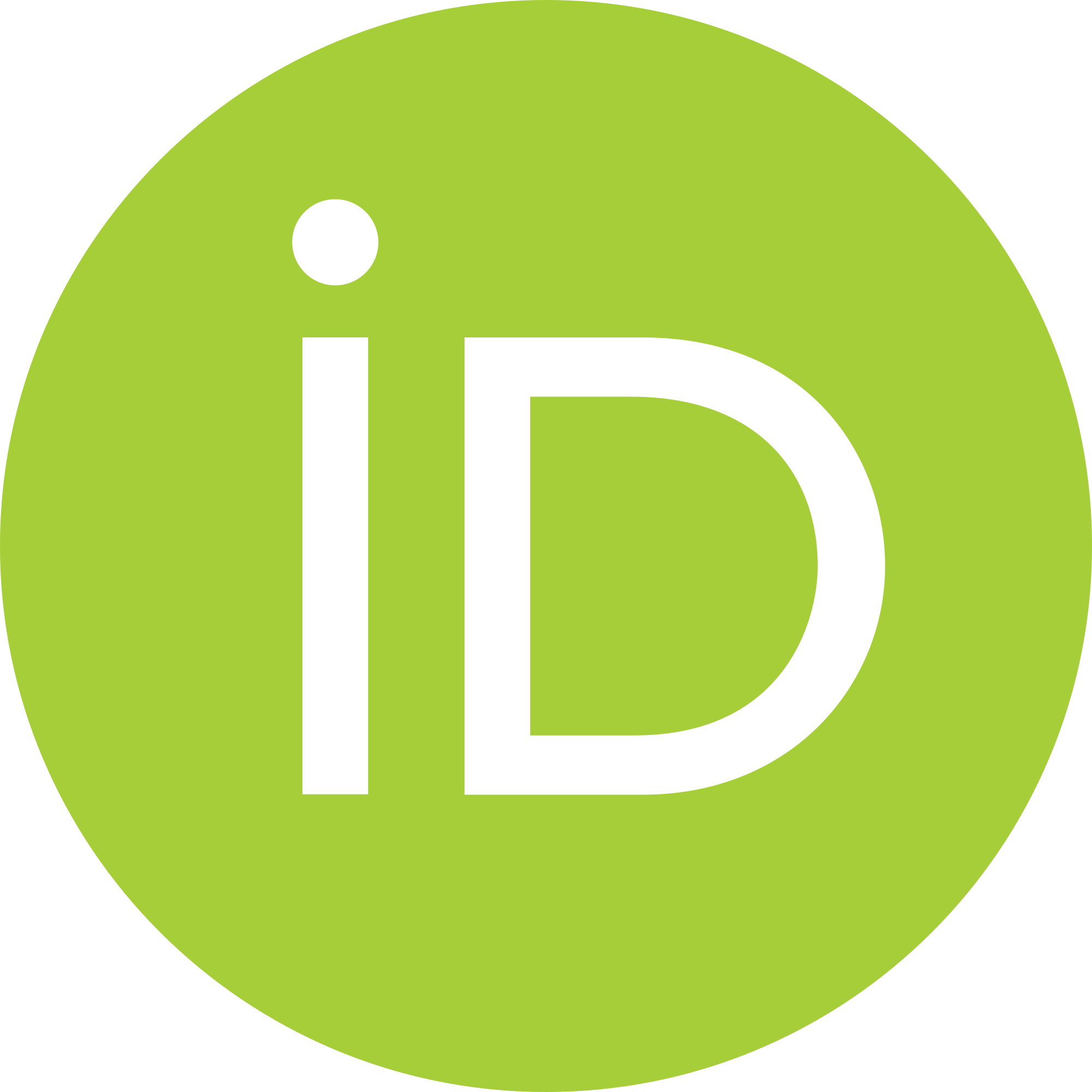}}}

\usepackage{verbatim}
\usepackage{hyperref}
\hypersetup{colorlinks=true,linkcolor=blue,citecolor=blue,urlcolor=blue}
\DeclareMathOperator{\Irr}{Irr}

\DeclareMathOperator{\Hom}{Hom }
\DeclareMathOperator{\End}{End }

\DeclareMathOperator{\Vect}{\normalfont\underline{Vect}}
\DeclareMathOperator{\CC}{\C \boxtimes \overline{\C}}
\DeclareMathOperator{\tr}{tr}
\DeclareMathOperator{\Yon}{\yen}
\DeclareMathOperator{\id}{id}
\DeclareMathOperator{\Blank}{--}

\DeclareMathOperator{\Obj}{Obj}
\DeclareMathOperator{\Rep}{Rep}

\DeclareMathOperator{\tid}{\normalfont \textbf{1}}
\DeclareMathOperator{\fld}{\mathbb{K}}

\newcommand{\B}{\mathcal{B}}
\newcommand{\C}{\mathcal{C}}
\newcommand{\F}{\mathcal{F}}
\newcommand{\G}{\mathcal{G}}
\newcommand{\I}{\mathcal{I}}

\newcommand{\cH}{\mathcal{H}}
\newcommand{\RB}{\mathcal{RB}}
\newcommand{\RTA}{\mathcal{RTA}}
\newcommand{\RTC}{\mathcal{RTC}}
\newcommand{\TA}{\mathcal{TA}}
\newcommand{\TC}{\mathcal{TC}}

\newcommand{\<}{\langle}
\renewcommand{\>}{\rangle}

\newcommand{\idem}{\varepsilon}
\newcommand{\dual}{^*}
\newcommand{\Fun}{\operatorname{Fun}}
\newcommand{\refequals}[1][]{\stackrel{\text{#1}}{=} }
\newcommand{\op}{^{\mathrm{op}}}
\newcommand{\lra}[1]{\stackrel{#1}{\longrightarrow}}
\newcommand{\defeq}{\vcentcolon=}

\theoremstyle{plain}
\newtheorem{THM}{Theorem}[section]
\newtheorem{PROP}[THM]{Proposition}
\newtheorem{LEMMA}[THM]{Lemma}

\newtheorem{COR}[THM]{Corollary}
\theoremstyle{definition}
\newtheorem{REM}[THM]{Remark}

\newtheorem{DEF}[THM]{Definition}
\newtheorem*{DEF*}{Definition}

\numberwithin{equation}{section}

\title{Decomposing the Tube Category}

\begin{document}

\maketitle

	\begin{center}

	LEONARD HARDIMAN\orcidicon{0000-0003-1986-6704} and ALASTAIR KING\orcidicon{0000-0002-8399-2345} \vspace{0.3em}

	\textit{\footnotesize  Department of Mathematical Sciences, University of Bath, Bath, BA2 7AY, United Kingdom \\
	emails:  \url{leonard.p.a.hardiman@bath.edu}, \url{a.d.king@bath.ac.uk} }

	\end{center}

\begin{abstract}
The tube category of a modular tensor category is a variant of the tube algebra, first introduced by Ocneanu.
As a category, it can be decomposed in two different, but related, senses.
Firstly, via the Yoneda embedding, the Hom spaces decompose into summands factoring though irreducible functors,
in a manner analogous to decomposing an algebra as a sum of matrix algebras.
We describe these summands.
Secondly, under the Yoneda embedding, each object decomposes into irreducibles, which correspond to primitive idempotents in the category itself.
We identify these idempotents.
We make extensive use of diagram calculus in the description and proof of these decompositions.
\end{abstract}

\vspace{0.1cm}


\thispagestyle{empty}

\section{Introduction}

The \emph{tube algebra} of a monoidal category $ \C $ was first introduced by Ocneanu \cite{Ocneanu1993} in the realm of operator algebra theory. Connections between the tube algebra and the Drinfeld centre $ Z(\C) $ have been the subject of much research, culminating in a result of Popa, Shlyakhtenko and Vaes \cite[Proposition 3.14]{Vaes} that the category of representations of the tube algebra is equivalent to $ Z(\C) $. In the case when $ \C $ is a modular tensor category (MTC) $ Z(\C) $ is equivalent to the category $ \CC $ where $ \boxtimes $ denotes the Deligne tensor product and $ \overline{\C} $ is obtained by equipping $ \C $ with the opposite braiding (see \cite{MR1990929} or \cite[Proposition 8.20.12]{EtingofBook}). \\

Here we will focus instead on the \emph{tube category} $ \TC $, which is a multiobject version of the tube algebra and is Morita equivalent to it. In other words, the category of representations $ \RTC \defeq \Fun(\TC^{\op}, \Vect) $ is equivalent to that of the tube algebra, where $ \Vect $ denotes the category of finite dimensional vector spaces. Thus, by the results mentioned above, when $ \C $ is an MTC we have
	\begin{align*}
	\RTC \cong Z(\C) \cong \CC.
	\end{align*}
As $ \{ I \boxtimes J \}_{I,J \in \Irr(\C)} $ forms a complete set of simple objects in $ \CC $, this equivalence gives a complete set $ \{F_{IJ}\}_{I,J \in \Irr(\C)} $ of  irreducible functors in $ \RTC $. Using the Yoneda embedding, which maps an object $ X $ in $ \TC $ to $ X^\sharp = \Hom_{\TC}(\Blank,X) $ in $ \RTC $, one obtains the decomposition
	\begin{align} \label{EQ:DECOMPOSE}
	\Hom_{\TC}(X,Y) &= \Hom_{\RTC}(X^\sharp,Y^\sharp) \nonumber \\
	&= \bigoplus\limits_{I,J} \Hom_{\RTC} (X^\sharp, F_{IJ}) \otimes \Hom_{\RTC} (F_{IJ},Y^\sharp).
	\end{align}
For fixed $ I,J \in \Irr(\C) $, $ \Hom_{\RTC} (X^\sharp, F_{IJ})  $ can be identified with $ F_{IJ}(X) $ via the canonical Yoneda map and $ \Hom_{\RTC}(F_{IJ},Y^\sharp) $ can first be identified with $ \Hom_{\RTC}(Y^\sharp,F_{IJ})\dual$ via the perfect pairing given by composition, then identified with $ F_{IJ}(Y)\dual $ via Yoneda again. Putting all of this together the $ (I,J) $ summand of \eqref{EQ:DECOMPOSE} can be identified with $ F_{IJ}(Y)\dual \otimes F_{IJ}(X) $, giving rise to a natural injection
	\begin{align*}
	\lambda_{YX}^{IJ} \colon F_{IJ}(Y)\dual \otimes F_{IJ}(X) \to \Hom_{\TC}(X,Y)
	\end{align*}
which, by naturality in $ X $, may also be thought of as a map
	\begin{align*}
	\lambda_{Y}^{IJ} : \F_{IJ}(Y)^* \to \Hom_{\RTC}(\F_{IJ},Y^\sharp).
	\end{align*}

The main goal of this paper is to give a graphical description of  $ \lambda_{YX}^{IJ} $ and thus $ \lambda_{Y}^{IJ} $. Composition is easily described with respect to this map (see Proposition~\ref{PROP:compo}) and this allows us to identify the primitive idempotents in $ \End_{\TC}(X) $ (see Corollary \ref{COR:PRIMITIVE}) which correspond to the irreducible summands of $ X^\sharp $ in $ \RTC $. These idempotents may be thought of as categorical analogues of \emph{Ocneanu projections} \cite{Ocneanu1993}, \cite[Section 2]{Evans1998}.\\

The structure of this paper is as follows. Section~\ref{SEC:PRELIMSONYON} starts by developing some general formalism, in particular Lemma~\ref{lem:compo} implies that $ \lambda_{Y}^{IJ} $ is characterised by being the unique opposite (see Definition~\ref{DEF:OPPOSITE}) of the canonical Yoneda map
	\begin{align*}
	\mu_Y^{IJ} \colon F_{IJ}(Y) \lra{\cong} \Hom_{\RTC}(Y^\sharp, F_{IJ}).
	\end{align*}
In Section~\ref{SEC:GRAPHICALCALCULUS} we build on the graphical calculus of MTC's, culminating in two main results: Lemma~\ref{LEMMA:DUAL_DECOMPOSE} and Proposition~\ref{PROP:TWISTED_S}.  Section~\ref{SEC:INTROTC} then provides an introduction to the tube category and finally Section~\ref{SECTCRTC} gives a graphical candidate for  $ \lambda_{Y}^{IJ} $ and proves that it is opposite to $ \mu^{IJ}_Y $, using the results of Section~\ref{SEC:GRAPHICALCALCULUS}. \\ \\
\textbf{Acknowledgement.} This work grew out of multiple stimulating discussions that took place during the authors' participation in the workshop ``Structure of operator algebras: subfactors and fusion categories'' held at the Isaac Newton Institute for Mathematical Sciences. We would like to thank the organisers for their hospitality and the participants for their helpful engagement. We would also like to extend particular thanks to Corey~Jones for enlightening conversation.
\section{Preliminaries on Yoneda in a $ \fld $-linear Category} \label{SEC:PRELIMSONYON}

We start by recording some basic results, which arise out of the Yoneda Lemma (see \cite[Section III.2]{MR1712872}) and that we use throughout. Let $ \fld $ be a field, let $\B$ be any $\fld$-linear category and let $\RB = \Rep\B = \Fun(\B\op,\Vect)$ be the category of $\fld$-linear contravariant functors from $ \B $ to $ \Vect $. We consider the Yoneda embedding
	\begin{align}\label{eq:Yonemb}
	\begin{split}
	\Yon \colon \B &\to\RB \\
	X &\mapsto X^\sharp
	\end{split}
	\end{align}
where $X^\sharp = \Hom_{\B}(\Blank,X)$. For any $F$ in $\RB$ and any $X$ in $\B$,
the Yoneda Lemma gives a natural isomorphism
	\begin{align}\label{eq:yon1mu}
	 \mu_X \colon F(X) \lra{\cong} \Hom_\RB(X^\sharp, F),
	\end{align}
where, for $\alpha\in F(X)$ and any $\eta\in X^\sharp(Z)$, i.e. $\eta\colon Z\to X$, we have
	\begin{align*}
	\mu_{XZ}(\alpha,\eta)\defeq\mu_X(\alpha)_Z(\eta)=F(\eta)(\alpha).
	\end{align*}
In other words, for any $Z$ in $\B$, we have a (bi)linear map
	\begin{align}\label{eq:yon2mu}
	\begin{split}
	\mu_{XZ}\colon F(X)\otimes \Hom_\B(Z,X) &\to F(Z)\\
	\alpha\otimes \eta &\mapsto F(\eta)(\alpha).
	\end{split}
	\end{align}
The inverse $\mu_X^{-1}\colon \Hom_\RB(X^\sharp, F) \to F(X)$ is given by $\Phi\mapsto \Phi_X(\id)$.
To see this, observe that
	\begin{align}\label{eq:mu-inv}
	\mu_X(\alpha)_X(\id)=F(\id)(\alpha)=\alpha,
	\end{align}
while the naturality of $\Phi$ implies that, for any $\eta\colon Z\to X$,
	\begin{align*}
	\mu_X(\Phi_X(\id))_Z(\eta)=F(\eta)(\Phi_X(\id))=\Phi_Z(\eta^*(\id))=\Phi_Z(\eta),
	\end{align*}
that is,  $\mu_X(\Phi_X(\id))=\Phi$.
	\begin{DEF} \label{DEF:OPPOSITE}
	For $ \mu_X $ as in \eqref{eq:yon1mu}, a map
		\begin{align}\label{eq:yon1lam}
		\lambda_X\colon F(X)\dual \lra{} \Hom_\RB(F,X^\sharp)
		\end{align}
	is said to be an \emph{opposite of $ \mu_X $} if it is formally dual to $\mu_X^{-1}$ in the sense that
	\begin{equation}\label{eq:dual}
	\mu_X(\alpha)\circ \lambda_X(\beta) = \<\beta,\alpha\> \id_F,
	\end{equation}
	where $\<\beta,\alpha\>=\beta(\alpha)$ is the natural duality pairing and $\id_F\in\End_\RB(F)$ is the identity natural transformation.
	\end{DEF}
	\begin{REM} \label{REM:SCHUR}
	If $F$ is simple and Schurian, i.e.  $\End_\RB(F)=\fld$, and $\RB$ is semisimple, then $\Hom_\RB(F,X^\sharp)$ and $\Hom_\RB(X^\sharp,F)$ are dual and the pairing is given by composition into $\End_\RB(F)$ (c.f. Proposition~\ref{PROP:SIMPLE_PERFECT}). In that case, \eqref{eq:dual} uniquely characterises $\lambda_X(\beta)$:
it is precisely the dual of $\mu_X^{-1}$.
	\end{REM}
We now suppose we have such a map $ \lambda_X $. As for $\mu$, we can write a bilinear version of $\lambda$ as
	\begin{align}\label{eq:yon2lam}
	\begin{split}
	\lambda_{XZ}\colon F(X)\dual \otimes F(Z) &\to \Hom_\B(Z,X)\\
	\beta\otimes \gamma &\mapsto \lambda_X(\beta)_Z(\gamma).
	\end{split}
	\end{align}
We will also denote the result of this map by $\lambda_{XZ}(\beta,\gamma)$. We consider the following chain of maps
\begin{align}
\label{eq:chain}
\begin{split}
F(Y)\dual \otimes &F(X)  \lra{\lambda_Y\otimes\mu_X} \Hom_\RB(F,Y^\sharp)\otimes\Hom_\RB(X^\sharp, F) \\
&\lra{\circ} \Hom_\RB(X^\sharp,Y^\sharp) \lra{} \Hom_\B(X,Y)
\end{split}
\end{align}
where the last map is the inverse of the isomorphism \eqref{eq:yon1mu} with $F=Y^\sharp$,
that is, the inverse of
$\Yon\colon \Hom_\B(X,Y)\to \Hom_\RB(X^\sharp,Y^\sharp)\colon \phi\mapsto\phi^\sharp$,
whose existence certifies that the Yoneda embedding is fully faithful.
	\begin{LEMMA}\label{lem:compo}
	The composition of this chain is $\lambda_{YX}$. In other words, in $\RB$ we have, for any $\beta\in F(Y)\dual$ and $\alpha\in F(X)$,
	\begin{align}\label{eq:compo}
	\lambda_{YX}(\beta,\alpha)^\sharp = \lambda_Y(\beta)\circ \mu_X(\alpha),
	\end{align}
	as a composition of maps $X^\sharp\to F\to Y^\sharp$.
	\proof
	In the case $\Phi=\lambda_Y(\beta)\circ \mu_X(\alpha)$, we compute
	\begin{align*}
	\Phi_X(\id) = \lambda_Y(\beta)_X\bigl(  \mu_X(\alpha)_X(\id)\bigr)
	= \lambda_Y(\beta)_X(\alpha)
	\end{align*}
	by \eqref{eq:mu-inv}, as required. \endproof
	\end{LEMMA}
	\begin{PROP}\label{PROP:compo}
	As a composition of maps $X\to Y\to Z$ in $\B$, we have
	\begin{align*}
	\lambda_{ZY}(\delta,\gamma) \circ \lambda_{YX}(\beta, \alpha)
	= \<\beta,\gamma\> \lambda_{ZX}(\delta,\alpha)
	\end{align*}
	\proof 	Under the Yoneda embedding and using Lemma~\ref{lem:compo}, the equation is the same,
	as a composite of maps $X^\sharp\to Y^\sharp\to Z^\sharp$ in $\RB$, as
	\begin{align*}
	\lambda_Z(\delta)\circ \mu_Y(\gamma)\circ
	\lambda_Y(\beta)\circ \mu_X(\alpha)
	= \<\beta,\gamma\> \lambda_Z(\delta)\circ \mu_X(\alpha),
	\end{align*}
	which follows immediately by applying the duality \eqref{eq:dual}. \endproof
	\end{PROP}

For any idempotent $\idem\in\End_\B(X)$ there is a subfunctor $(X,\idem)^\sharp\leq X^\sharp$ in $\RB$ which is the image of the idempotent $\idem^\sharp\in \End_\RB(X^\sharp)$. Indeed, $(X,\idem)^\sharp$ is a summand of $X^\sharp$. This image exists because $\RB$ is an abelian category (see, for example, \cite[Section VIII.3]{MR1712872}), so idempotent complete, even if $\B$ may not be. Concretely, $(X,\idem)^\sharp(Y)$ is the image of $\idem^\sharp_Y=\idem_*$, which is an idempotent endomorphism of $X^\sharp(Y)=\Hom_\B(Y,X)$. The naturality of $\idem^\sharp$, i.e. the fact that $\idem_*$ commutes with $\phi^*$ for any $\phi\colon Z\to Y$, makes $(X,\idem)^\sharp$ a functor.\\

For any object $X$ in $ \B$ together with $\beta\in F(X)\dual$ and $\alpha\in F(X)$
such that $\<\beta,\alpha\>=1$, Corollary~\ref{PROP:compo} implies that $ \psi = \lambda_{XX}(\beta,\alpha)$ is an idempotent in $\End_\B(X)$.
	\begin{COR} \label{COR:isoidem}
	We have $ (X,\psi)^\sharp \cong  F$.
	\proof  By \eqref{eq:compo}, we have that
	$\psi^\sharp= \lambda_X(\beta)\circ\mu_X(\alpha)$ and, by \eqref{eq:dual}, that
	$\mu_X(\alpha)\circ \lambda_X(\beta)=\<\beta,\alpha\> \id_F=\id_F$, by assumption.
	Hence, $(X,\psi)^\sharp\cong F$ in a way that identifies
	$\mu_X(\alpha)\colon X^\sharp\to F$ and $\lambda_X(\beta)\colon F\to X^\sharp$
	with projection and inclusion of the summand. \endproof
	\end{COR}
\section{Graphical Calculus in Modular Tensor Categories} \label{SEC:GRAPHICALCALCULUS}

Throughout this section we make extensive use of the graphical calculus described in \cite[Section 2.3]{GraphCal}. Our conventions are as follows: composition is read vertically top to bottom and tensor products are read horizontally left to right. The results in this section are well-known to experts in the field however we include their statements and proofs for the sake of completeness. \\

Let $ \fld $ be an algebraically closed field. From now on $ \C $ is always assumed to be a \emph{modular tensor category} over $ \fld $ (see \cite[Definition 3.1.1]{GraphCal} for a precise definition). We recall this implies that $ \C $ is semisimple of finite type and that the tensor identity, denoted $ \tid $, is a simple object. We use $ \Irr(\C) $ to denote a complete set of simple objects in $ \C $. We have $ \End_{\C}(S) = \fld $ for all $ S \in \Irr(\C) $ by the algebraic closure of $ \mathbb{K} $ and Schur's Lemma. Unless otherwise specified a sum over a variable object in $ \C $ ranges over $ \Irr(\C) $.
	\begin{PROP} \label{PROP:SIMPLE_PERFECT}
	Let $ R $ be in $ \Irr(\C) $ and let $ X $ be in $ \C $. The pairing
	\begin{align*}
	\Hom_{\C}(R,X) \otimes \Hom_{\C}(X,R) &\to \fld \\
	f \otimes g &\mapsto g \circ f
	\end{align*}
	is perfect.
	\proof By semisimplicity we have $ X = \bigoplus_{i \in \I} X_i $ where the $ X_i $ are simple objects and $ \I $ is an indexing set. We consider the subset $ J \defeq \{ i \in I \mid X_i \cong R \} \subset I $. Then we have
		\begin{align*}
		\Hom_{\C}(R,X) \cong \fld^J, \ \Hom_{\C}(X,R) \cong \fld^J
		\end{align*}
	and composition is given by the standard pairing, which is perfect.\endproof
	\end{PROP}
	\begin{DEF}
	Let $ R $ be in $ \Irr(\C) $ and let $ X $ be in $ \C $. For every basis $ \{b\} \subset \Hom_{\C}(R,X) $ we use $ \{b^*\} $ to denote the dual basis of $ \Hom_{\C}(X,R) $ with respect to the perfect pairing given by Proposition~\ref{PROP:SIMPLE_PERFECT}.
	\end{DEF}
	\begin{LEMMA}\label{LEMMA:DECOMPOSE}
	Let $ X $ be in $ C $. We have
			\begin{align} \label{EQ:LEMMA_DECOMPOSE}
			\begin{array}{c}
			\begin{tikzpicture}[scale=0.15,every node/.style={inner sep=0,outer sep=-1}]
			\node (v8) at (0,6) {};
			\node (v10) at (0,-6) {};
			\node at (-2,2) {$X$};
			\draw[thick]  (v8) edge (v10);
			\end{tikzpicture}
			\end{array}
			= \sum\limits_{R,b}
			\begin{array}{c}
			\begin{tikzpicture}[scale=0.15,every node/.style={inner sep=0,outer sep=-1}]
			\node (v2) at (0,7.6) {};
			\node (v4) at (0,-7.8) {};
			\node [draw,outer sep=0,inner sep=2,minimum width=14,minimum height = 12] (v6) at (0,3.5) {$ b^* $};
			\node [draw,outer sep=0,inner sep=2,minimum width=14,minimum height = 12] (v9) at (0,-3.5) {$ b $};
			\node at (-1.75,6.5) {$X$};
			\node at (-1.75,-6.75) {$X$};
			\node at (1.75,0) {$R$};
			\draw[thick]  (v6) edge (v9);
			\draw[thick]  (v2) edge (v6);
			\draw[thick]  (v9) edge (v4);
			\end{tikzpicture}
			\end{array}
			\end{align}
		where $ b $ ranges over a basis of $ \Hom_{\C}(R,X) $.
	\proof
By semisimplicity, we have a natural identification $ \bigoplus_R \Hom_{\C}(R,X) \otimes R = X  $.
Using the bases $ \{ b \} $, for each $ \Hom_{\C}(R,X) $, we get an isomorphism
$ f \colon \bigoplus_{R,b} R \to X $, given diagrammatically by
		\begin{align*}
	     f_{R,b} =
		\begin{array}{c}
		\begin{tikzpicture}[scale=0.15,every node/.style={inner sep=0,outer sep=-1}]
		\node [draw,outer sep=0,inner sep=2,minimum width=14,minimum height = 12] (v9) at (0,-3) {$ b $};
		\node at (1.5,-6.8) {$X$};
		\node at (1.5,0.5) {$R$};
		\draw[thick]  (0,2) edge (v9);
		\draw[thick]  (v9) edge (0,-8);
		\end{tikzpicture}
		\end{array}.
		\end{align*}
	On the other hand, the map from $ X $ to $ \bigoplus_{R,b} R $ given by
		\begin{align*}
	     g_{R,b} =
		\begin{array}{c}
		\begin{tikzpicture}[scale=0.15,every node/.style={inner sep=0,outer sep=-1}]
		\node [draw,outer sep=0,inner sep=2,minimum size=10] (v9) at (0,-3) {$ b\dual $};
		\node at (1.5,-6.8) {$R$};
		\node at (1.5,0.5) {$X$};
		\draw[thick]  (0,2) edge (v9);
		\draw[thick]  (v9) edge (0,-8);
		\end{tikzpicture}
		\end{array}
		\end{align*}
	satisfies $ g \circ f = \id $ and so $ g $ is the inverse to $ f $. The right-hand side of \eqref{EQ:LEMMA_DECOMPOSE} is simply $ f \circ g $ and hence equal to the identity, as required.
	\endproof
	\end{LEMMA}
	\begin{REM} \label{REM:SPECIALISED_LEMMA_DECOMPOSE}
		We will often use the following instance of Lemma~\ref{LEMMA:DECOMPOSE}. For $ S,T \in \Irr(\C) $ we have
		\begin{align*}
		\begin{array}{c}
		\begin{tikzpicture}[scale=0.15,every node/.style={inner sep=0,outer sep=-1}]
		\node (v5) at (-2,6) {};
		\node (v7) at (-2,-6) {};
		\node (v8) at (2,6) {};
		\node (v10) at (2,-6) {};
		\node at (-3.5,2) {$S$};
		\node at (3.5,2) {$T$};
		\draw[thick]  (v5) edge (v7);
		\draw[thick]  (v8) edge (v10);
		\end{tikzpicture}
		\end{array}
		= \sum\limits_{R,b}
		\begin{array}{c}
		\begin{tikzpicture}[scale=0.15,every node/.style={inner sep=0,outer sep=-1}]
		\node (v2) at (4,6.5) {};
		\node (v1) at (-4,6.5) {};
		\node (v4) at (4,-6.5) {};
		\node (v3) at (-4,-6.5) {};
		\node (v5) at (-2.1,6.5) {};
		\node (v7) at (-0.5,4) {};
		\node (v8) at (2.1,6.5) {};
		\node (v10) at (0.5,4) {};
		\node (v11) at (-2.1,-6.5) {};
		\node (v13) at (-0.5,-4) {};
		\node (v14) at (2.1,-6.5) {};
		\node (v16) at (0.5,-4) {};
		\node at (-3.5,5.5) {$S$};
		\node at (3.5,5.5) {$T$};
		\node at (-3.5,-5.5) {$S$};
		\node at (3.5,-5.5) {$T$};
		\node at (1.5,0) {$R$};
		\draw[thick]  plot[smooth, tension=.7] coordinates {(v7) (-1,5) (-1.9,5.5) (v5)};
		\draw[thick]  plot[smooth, tension=.7] coordinates {(v10) (0.9,5) (1.9,5.5) (v8)};
		\draw[thick]  plot[smooth, tension=.7] coordinates {(v13) (-0.9,-5) (-1.9,-5.5) (v11)};
		\draw[thick]  plot[smooth, tension=.7] coordinates {(v16) (0.9,-5) (1.9,-5.5) (v14)};
		\node [draw,outer sep=0,inner sep=2,minimum size=10,fill=white] (v6) at (0,3) {$ b^* $};
		\node [draw,outer sep=0,inner sep=2,minimum width=14,minimum height = 12,fill=white] (v9) at (0,-3) {$ b $};
		\draw[thick]  (v6) edge (v9);
		\end{tikzpicture}
		\end{array}
		\end{align*}
	where $ b $ ranges over a basis of $ \Hom_{\C}(R,ST) $.
	\end{REM}
We recall that modular tensor categories are \emph{rigid} i.e. any object $ X $ in $ \C $ admits a left dual and a right dual. We do not distingish between these two objects as they are identified by the \emph{pivotal structure} of $ \C $, we denote them both by $ X^\vee $. The corresponding structural maps are diagrammatically represented as
	$$
	\begin{array}{c}
	\begin{tikzpicture}
	\node at (-0.5,-0.1) {$X$};
	\node at (0.5,-0.1) {$X^\vee$};
	\node (v1) at (-0.5,0) {};
	\node (v2) at (0.5,0) {};
	\draw[thick] (v1) to[out=90,in=90] (v2);
	\end{tikzpicture}
	\end{array} \text{and}
	\begin{array}{c}
	\begin{tikzpicture}
	\node at (-0.5,0.1) {$X$};
	\node at (0.5,0.1) {$X^\vee$};
	\node (v1) at (-0.5,0) {};
	\node (v2) at (0.5,0) {};
	\draw[thick] (v1) to[out=-90,in=-90] (v2);
	\end{tikzpicture}.
	\end{array} $$
	\begin{DEF}
	Let $ f $ be in $ \End_{\C}(X) $ for some object $ X $ in $ \C $. The \emph{trace} of $ f $ is defined by
		\begin{align*}
		\tr:
		\begin{array}{c}
		\begin{tikzpicture}[scale=0.15,every node/.style={inner sep=0,outer sep=-1}]
		\node (v2) at (4,6) {};
		\node (v1) at (-4,6) {};
		\node (v4) at (4,-6) {};
		\node (v3) at (-4,-6) {};
		\node (v5) at (0,6) {};
		\node (v7) at (0,-6) {};
		\node at (2,5) {$X$};
		\node at (2,-5) {$X$};
		\node [draw,outer sep=0,inner sep=3,minimum size=10] (v6) at (0,0) {$ f $};
		\draw[thick]  (v6) edge (v5);
		\draw[thick]  (v7) edge (v6);
		\end{tikzpicture}
		\end{array}
		\longmapsto \ \
		\begin{array}{c}
		\begin{tikzpicture}[scale=0.15,every node/.style={inner sep=0,outer sep=-1}]
		\node at (4.5,0) {$ X^\vee $};
		\node [draw,outer sep=0,inner sep=3,minimum size=10] (v6) at (-1,0) {$ f $};
		\node (v8) at (-1,3) {};
		\node (v10) at (2,3) {};
		\node (v9) at (-1,-3) {};
		\node (v11) at (2,-3) {};
		\draw[thick]  (v8) edge (v6);
		\draw[thick]  (v9) edge (v6);
		\draw[thick] (v8) to[out=90,in=90] (v10);
		\draw[thick]  (v10) edge (v11);
		\draw[thick] (v9) to[out=-90,in=-90] (v11);
		\end{tikzpicture}
		\end{array} \in \End_{\C}(\tid) = \fld.
		\end{align*}
	We recall that modular tensor categories are \emph{spherical} i.e.
		\begin{align*}
		\begin{array}{c}
		\begin{tikzpicture}[scale=0.15,every node/.style={inner sep=0,outer sep=-1},xscale=-1]
		\node at (4.5,0) {$ X^\vee $};
		\node [draw,outer sep=0,inner sep=3,minimum size=10] (v6) at (-1,0) {$ f $};
		\node (v8) at (-1,3) {};
		\node (v10) at (2,3) {};
		\node (v9) at (-1,-3) {};
		\node (v11) at (2,-3) {};
		\draw[thick]  (v8) edge (v6);
		\draw[thick]  (v9) edge (v6);
		\draw[thick] (v8) to[out=90,in=90] (v10);
		\draw[thick]  (v10) edge (v11);
		\draw[thick] (v9) to[out=-90,in=-90] (v11);
		\end{tikzpicture}
		\end{array} =
		\begin{array}{c}
		\begin{tikzpicture}[scale=0.15,every node/.style={inner sep=0,outer sep=-1}]
		\node at (4.5,0) {$ X^\vee $};
		\node [draw,outer sep=0,inner sep=3,minimum size=10] (v6) at (-1,0) {$ f $};
		\node (v8) at (-1,3) {};
		\node (v10) at (2,3) {};
		\node (v9) at (-1,-3) {};
		\node (v11) at (2,-3) {};
		\draw[thick]  (v8) edge (v6);
		\draw[thick]  (v9) edge (v6);
		\draw[thick] (v8) to[out=90,in=90] (v10);
		\draw[thick]  (v10) edge (v11);
		\draw[thick] (v9) to[out=-90,in=-90] (v11);
		\end{tikzpicture}
		\end{array}
		\end{align*}
	 and that taking the trace of a composition is commutative i.e.
		\begin{align*}
		\begin{array}{c}
		\begin{tikzpicture}[scale=0.15,every node/.style={inner sep=0,outer sep=-1}]
		\node [draw,outer sep=0,inner sep=2.4,minimum height=15,minimum width=11] (v9) at (0,2) {$ g $};
		\node [draw,outer sep=0,inner sep=2.2,minimum height=15,minimum width=11] (v90) at (0,9.5) {$ f $};
		\node at (-6.25,5.75) {$X^\vee$};
		\node at (-1.5,5.75) {$Y$};
		\draw[thick]  (v90) edge (v9);
		\node (v1) at (-4,0) {};
		\node (v2) at (-4,11.5) {};
		\draw [thick]  (v9) to[out=-90,in=-90] (v1);
		\draw [thick] (v90) to[out=90,in=90] (v2);
		\draw [thick] (v2) edge (v1);
		\end{tikzpicture}
		\end{array} =
		\begin{array}{c}
		\begin{tikzpicture}[scale=0.15,every node/.style={inner sep=0,outer sep=-1}]
		\node [draw,outer sep=0,inner sep=2.4,minimum height=15,minimum width=11] (v9) at (0,2) {$ f $};
		\node [draw,outer sep=0,inner sep=2.2,minimum height=15,minimum width=11] (v90) at (0,9.5) {$ g $};
		\node at (-2.5,5.75) {$Y$};
		\node at (2.5,5.75) {$X^\vee$};
		\draw[thick]  (v90) edge (v9);
		\node (v1) at (-4,0) {};
		\node (v2) at (-4,11.5) {};
		\draw [thick]  (v9) to[out=-90,in=-90] (v1);
		\draw [thick] (v90) to[out=90,in=90] (v2);
		\draw [thick] (v2) edge (v1);
		\end{tikzpicture}
		\end{array}
		\end{align*}
	for all $ f \in \Hom_{\C}(X,Y) $ and $ g \in \Hom_{\C}(Y,X) $.
	\end{DEF}
	\begin{DEF}
	Let $ X $ be an object in $ \C $. The \emph{dimension} of $ X $ is defined by
		\begin{align*}
		d(X) = \tr(\id_X).
		\end{align*}
	\end{DEF}
	\begin{DEF}
	The \emph{dimension} of $ \C $ is defined by
		\begin{align*}
		d(\C) = \sum\limits_S d(S)^2.
		\end{align*}
	\end{DEF}
	\begin{REM}
	From \cite[Proposition 4.8.4.]{EtingofBook} we have that $ S \in \Irr(\C) $ implies $ d(S) \neq 0 $. From \cite[Proposition 8.20.16.]{EtingofBook} we have $ d(\C) \neq 0 $ (which requires $ \mathbb{K} $ to be algebraically closed). We will use both facts repeatedly
throughout.
	\end{REM}
	\begin{PROP} \label{PROP:NONSIMPLE_PERFECT}
	The pairing
	\begin{align*}
	\Hom_{\C}(X,Y) \otimes \Hom_{\C}(Y,X) &\to \fld \\
	f \otimes g &\mapsto \tr(g \circ f)
	\end{align*}
	is perfect. \proof
	We have
		\begin{align*}
		\tr(g \circ f) &=
		\begin{array}{c}
		\begin{tikzpicture}[scale=0.15,every node/.style={inner sep=0,outer sep=-1}]
		\node [draw,outer sep=0,inner sep=2.4,minimum height=15] (v9) at (-0.25,3)  {$ g $};
		\node [draw,outer sep=0,inner sep=2.2,minimum size=10] (v90) at (-0.25,9) {$ f $};
		\node at (-6.8,5.9) {$X^\vee$};
		\node at (3,5.9) {$Y$};
		\draw[thick]  (v90) edge (v9);
		\node (v1) at (-4,1.1) {};
		\node (v2) at (-4,10.75) {};
		\draw [thick]  (v9) to[out=-90,in=-90] (v1);
		\draw [thick] (v90) to[out=90,in=90] (v2);
		\draw [thick] (v2) edge (v1);
		\end{tikzpicture}
		\end{array}~\refequals[Lem.~\ref{LEMMA:DECOMPOSE}] \sum\limits_{S,b}
		\begin{array}{c}
		\begin{tikzpicture}[scale=0.15,every node/.style={inner sep=0,outer sep=-1}]
		\node (v6) at (0,21.4) {};
		\node (v3) at (-4,21.4) {};
		\node (v5) at (0,1.1) {};
		\node (v4) at (-4,1.1) {};
		\node [draw,outer sep=0,inner sep=1,minimum height=15, minimum width=12] (v9) at (0,3) {$ g $};
		\node [draw,outer sep=0,inner sep=1,minimum height=15, minimum width=12] (v90) at (0,9) {$ f $};
		\node [draw,outer sep=0,inner sep=1,minimum height=15, minimum width=12] (v1) at (0,14.5) {$b$};
		\node [draw,outer sep=0,inner sep=1,minimum height=15, minimum width=12] (v2) at (0,19.5) {$b^*$};
		\node at (3.2,11.75) {$X$};
		\node at (-6.5,11.684) {$X^\vee$};
		\node at (3.2,17) {$S$};
		\node at (3.2,5.842) {$Y$};
		\draw[thick]  (v90) edge (v9);
		\draw [thick] (v90) edge (v1);
		\draw [thick] (v2) edge (v1);
		\draw [thick] (v3) edge (v4);
		\draw [thick] (v5) to[out=-90,in=-90] (v4);
		\draw [thick] (v6) to[out=90,in=90] (v3);
		\end{tikzpicture}
		\end{array} \\
		&= \sum\limits_{S,b}
		\begin{array}{c}
		\begin{tikzpicture}[scale=0.15,every node/.style={inner sep=0,outer sep=-1}]
		\node (v6) at (0,21.9) {};
		\node (v3) at (-4,21.9) {};
		\node (v5) at (0,1.1) {};
		\node (v4) at (-4,1.1) {};
		\node [draw,outer sep=0,inner sep=1.5,minimum height=15, minimum width=12] (v9) at (0,3) {$ b^* $};
		\node [draw,outer sep=0,inner sep=1,minimum height=15, minimum width=12] (v90) at (0,9) {$ g $};
		\node [draw,outer sep=0,inner sep=1,minimum height=15, minimum width=12] (v1) at (0,14.5) {$f$};
		\node [draw,outer sep=0,inner sep=1.5,minimum height=15, minimum width=12] (v2) at (0,20) {$b$};
		\node at (3.2,11.7) {$Y$};
		\node at (-6.096,11.684) {$S^\vee$};
		\node at (3.2,17.25) {$X$};
		\node at (3.2,6) {$X$};
		\draw[thick]  (v90) edge (v9);
		\draw [thick] (v90) edge (v1);
		\draw [thick] (v2) edge (v1);
		\draw [thick] (v3) edge (v4);
		\draw [thick] (v5) to[out=-90,in=-90] (v4);
		\draw [thick] (v6) to[out=90,in=90] (v3);
		\end{tikzpicture}
		\end{array} = \sum\limits_{S,b} \langle f \circ b , b^* \circ g \rangle d(S)
		\end{align*}
	where $ \langle\Blank,\Blank\rangle $ denotes the perfect pairing defined in Proposition~\ref{PROP:SIMPLE_PERFECT}. Let $ f \in \Hom_{\C}(X,Y) $ be non zero. Then there exists $ S \in \Irr(\C) $ and $ b \in \Hom_{\C}(S,X) $ such that $ f \circ b \neq 0 $. As $ \langle\Blank,\Blank\rangle $ is perfect there exists $ g_S \in \Hom_{\C}(Y,S) $ such that $ \langle f \circ b , g_S \rangle \neq 0 $. As $ d(S) \neq 0 $ we have
		\begin{align*}
		\tr(f \circ b \circ g_S) =  \langle f \circ b , b \circ b^* \circ g_S \rangle d(S) = \langle f \circ b , g_S \rangle d(S) \neq 0.
		\end{align*} \endproof
	\end{PROP}
	\begin{LEMMA} \label{LEMMA:TWISTED_DUALS}
 For any $ X,Y \in \C $, any $ S,S' \in \Irr(\C) $, any $ i \in \Hom_{\C}(Y,XS) $ and any $ j \in \Hom_{\C}(XS',Y) $ we have
		\begin{align*}
		\begin{array}{c}
		\begin{tikzpicture}[scale=0.15,every node/.style={inner sep=0,outer sep=-1}]
		\node (v6) at (0,7) {};
		\node (v2) at (0,6.5) {};
		\node [draw,outer sep=0,inner sep=2.2,minimum height=13, minimum width=9] (v9) at (0,3) {$ i $};
		\node [draw,outer sep=0,inner sep=1.5,minimum height=13, minimum width=9] (v90) at (0,9) {$ j $};
		\node at (-5.5,5.95) {$X^\vee$};
		\node at (2,0) {$S$};
		\node at (2.2,12) {$S'$};
		\node at (2,5.95) {$Y$};
		\draw[thick]  (v6) edge (v9);
		\draw [thick] (v90) edge (v2);
		\draw [thick] plot[smooth, tension=.7] coordinates {(0.5441,13) (0.5441,10.6)};
		\draw [thick] plot[smooth, tension=.7] coordinates {(0.5441,1.5) (0.5441,-1.1)};
		\node (v5) at (-3,10.625) {};
		\node (v1) at (-3,1.375) {};
		\node (v3) at (-0.5,1.375) {};
		\draw [thick] (v3) to[out=-90,in=-90] (v1);
		\node (v4) at (-0.5,10.625) {};
		\draw [thick] (v4) to[out=90,in=90] (v5);
		\draw [thick] (v5) edge (v1);
		\end{tikzpicture}
		\end{array}
		= \delta_{S,S'} \frac{\tr(j \circ i)}{d(S)}
		\begin{array}{c}
		\begin{tikzpicture}[scale=0.15,every node/.style={inner sep=0,outer sep=-1}]
		\draw[thick] (0,7) -- (0,-7);
		\node at (2,4) {$ S $};
		\end{tikzpicture}
		\end{array}.
		\end{align*}
	We also have, for any $ k \in \Hom_{\C}(Y,SX) $ and $ l \in \Hom_{\C}(S'X,Y) $,
	\begin{align*}
	\begin{array}{c}
	\begin{tikzpicture}[scale=0.15,every node/.style={inner sep=0,outer sep=-1},yscale=1,xscale=-1]
	\node (v6) at (0,7) {};
	\node (v2) at (0,6.5) {};
	\node [draw,outer sep=0,inner sep=2.2,minimum height = 13,minimum width = 11] (v9) at (0,3) {$ k $};
	\node [draw,outer sep=0,inner sep=1.5,minimum height = 13,minimum width = 11] (v90) at (0,9) {$ l $};
	\node at (-5.5,5.875) {$X^\vee$};
	\node at (2,0) {$S$};
	\node at (2.2,12) {$S'$};
	\node at (2,5.875) {$Y$};
	\draw[thick]  (v6) edge (v9);
	\draw [thick] (v90) edge (v2);
	\draw [thick] plot[smooth, tension=.7] coordinates {(0.5441,13) (0.5441,10.6)};
	\draw [thick] plot[smooth, tension=.7] coordinates {(0.5441,1.5) (0.5441,-1.1)};
	\node (v5) at (-3,10.6) {};
	\node (v1) at (-3,1.375) {};
	\node (v3) at (-0.5,1.375) {};
	\draw [thick] (v3) to[out=-90,in=-90] (v1);
	\node (v4) at (-0.5,10.6) {};
	\draw [thick] (v4) to[out=90,in=90] (v5);
	\draw [thick] (v5) edge (v1);
	\end{tikzpicture}
	\end{array}
	= \delta_{S,S'}  \frac{\tr(l \circ k)}{d(S)}
	\begin{array}{c}
	\begin{tikzpicture}[scale=0.15,every node/.style={inner sep=0,outer sep=-1}]
	\draw[thick] (0,7) -- (0,-7);
	\node at (2,4) {$ S $};
	\end{tikzpicture}
	\end{array}.
	\end{align*}
	\proof
	As $ S $ and $ S' $ are simple objects we have
		\begin{align*}
		\begin{array}{c}
		\begin{tikzpicture}[scale=0.15,every node/.style={inner sep=0,outer sep=-1}]
		\node (v6) at (0,7) {};
		\node (v2) at (0,6.5) {};
		\node [draw,outer sep=0,inner sep=2.2,minimum height=13, minimum width=9] (v9) at (0,3) {$ i $};
		\node [draw,outer sep=0,inner sep=1.5,minimum height=13, minimum width=9] (v90) at (0,9) {$ j $};
		\node at (-5.5,5.95) {$X^\vee$};
		\node at (2,0) {$S$};
		\node at (2.2,12) {$S'$};
		\node at (2,5.95) {$Y$};
		\draw[thick]  (v6) edge (v9);
		\draw [thick] (v90) edge (v2);
		\draw [thick] plot[smooth, tension=.7] coordinates {(0.5441,13) (0.5441,10.6)};
		\draw [thick] plot[smooth, tension=.7] coordinates {(0.5441,1.5) (0.5441,-1.1)};
		\node (v5) at (-3,10.625) {};
		\node (v1) at (-3,1.375) {};
		\node (v3) at (-0.5,1.375) {};
		\draw [thick] (v3) to[out=-90,in=-90] (v1);
		\node (v4) at (-0.5,10.625) {};
		\draw [thick] (v4) to[out=90,in=90] (v5);
		\draw [thick] (v5) edge (v1);
		\end{tikzpicture}
		\end{array}
		= \delta_{S,S'} \lambda
		\begin{array}{c}
		\begin{tikzpicture}[scale=0.15,every node/.style={inner sep=0,outer sep=-1}]
		\draw[thick] (0,7) -- (0,-7);
		\node at (2,4) {$ S $};
		\end{tikzpicture}
		\end{array}.
		\end{align*}
	for some $ \lambda \in \fld $. To compute $ \lambda $ we suppose $ \delta_{S,S'}=1 $ and take the trace to obtain
		\begin{align*}
		\lambda d(S) =
		\begin{array}{c}
		\begin{tikzpicture}[scale=0.15,every node/.style={inner sep=0,outer sep=-1}]
		\node (v6) at (0,7) {};
		\node (v2) at (0,6.5) {};
		\node [draw,outer sep=0,inner sep=2.2,minimum height=13, minimum width=9] (v9) at (0,3) {$ i $};
		\node [draw,outer sep=0,inner sep=1.5,minimum height=13, minimum width=9] (v90) at (0,9) {$ j $};
		\node at (-5.5,5.875) {$X^\vee$};
		\node at (2,0) {$S$};
		\node at (2.2,12) {$S$};
		\node at (2,5.875) {$Y$};
		\node at (-10,6) {$ S^\vee $};
		\draw[thick]  (v6) edge (v9);
		\draw [thick] (v90) edge (v2);
		\node (v5) at (-3,10.625) {};
		\node (v1) at (-3,1.375) {};
		\node (v3) at (-0.5,1.375) {};
		\draw [thick] (v3) to[out=-90,in=-90] (v1);
		\node (v4) at (-0.5,10.625) {};
		\draw [thick] (v4) to[out=90,in=90] (v5);
		\draw [thick] (v5) edge (v1);
		\node (v7) at (-8,10.7) {};
		\node (v8) at (-8,1.4) {};
		\draw [thick]  (v7) edge (v8);
		\node (v11) at (0.5,1.4) {};
		\node (v10) at (0.5,10.7) {};
		\draw [thick]  (v10) to[out=90,in=90] (v7);
		\draw [thick] (v8) to[out=-90,in=-90] (v11);
		\end{tikzpicture}
		\end{array}
		=
		\begin{array}{c}
		\begin{tikzpicture}[scale=0.15,every node/.style={inner sep=0,outer sep=-1}]
		\node [draw,outer sep=0,inner sep=1.5,minimum height=13, minimum width=9] (v9) at (-0.25,3) {$ j $};
		\node [draw,outer sep=0,inner sep=2.2,minimum height=13, minimum width=9] (v90) at (-0.25,9) {$ i $};
		\node at (-6.8,6) {$Y^\vee$};
		\node at (2,6) {$S$};
		\node at (-2.5,6) {$X$};
		\draw [thick] plot[smooth, tension=.7] coordinates {(-1,4.6) (-1,7.5)};
		\draw [thick] plot[smooth, tension=.7] coordinates {(0.5,4.6) (0.5,7.5)};
		\node (v1) at (-4.5,1.5) {};
		\node (v2) at (-4.5,10.5) {};
		\draw [thick] (v9) to[out=-90,in=-90] (v1);
		\draw [thick] (v2) edge (v1);
		\draw [thick] (v90)  to[out=90,in=90] (v2);
		\end{tikzpicture}
		\end{array}
		= \tr(j \circ i).
		\end{align*}
	Therefore $ \lambda = \frac{\tr(j \circ i)}{d(S)} $. This proves the first part of the lemma, the second part is proved analogously.  \endproof
	\end{LEMMA}
	\begin{LEMMA} \label{LEMMA:DUAL_DECOMPOSE}
	Let $ X $ be in $ \C $ and $ S $ be in $ \Irr(\C) $. We have
	\begin{align*}
	d(S)  \
	\begin{array}{c}
	\begin{tikzpicture}[scale=0.15,every node/.style={inner sep=0,outer sep=-1}]
	\node (v5) at (-2,6) {};
	\node (v7) at (-2,-6) {};
	\node (v8) at (2,6) {};
	\node (v10) at (2,-6) {};
	\node at (-4,2) {$X$};
	\node at (3.5,2) {$S$};
	\draw[thick]  (v5) edge (v7);
	\draw[thick]  (v8) edge (v10);
	\end{tikzpicture}
	\end{array}
	= \sum\limits_{T,b} d(T)
	\begin{array}{c}
	\begin{tikzpicture}[scale=0.15,every node/.style={inner sep=0,outer sep=-1}]
	\node (v13) at (-0.5,1.5) {};
	\node (v130) at (-0.5,-2) {};
	\node (v6) at (0,7) {};
	\node (v2) at (0,-7.5) {};
	\node at (3,-6.5) {$S$};
	\node at (-5.25,-6.5) {$X$};
	\node at (3,6) {$S$};
	\node at (3,-0.25) {$T$};
	\node at (-5.25,6) {$X$};
	\draw [thick] plot[smooth, tension=.7] coordinates {(0.5,2) (0.5,-2.5)};
	\node (v1) at (-3,1.5) {};
	\node (v4) at (-3,-2) {};
	\node (v5) at (-3,-7.5) {};
	\node (v3) at (-3,7) {};
	\draw [thick] (v13) to[out=-90,in=-90] (v1);
	\draw [thick] (v3) edge (v1);
	\draw [thick] (v130) to[out=90,in=90] (v4);
	\draw [thick] (v5) edge (v4);
	\node (v8) at (-0.5,-3) {};
	\node (v7) at (-0.5,2.5) {};
	\draw [thick] (v13) edge (v7);
	\draw [thick] (v130) edge (v8);
	\node [draw,outer sep=0,inner sep=2,minimum width=14,fill=white] (v9) at (0,3) {$ b $};
	\node [draw,outer sep=0,inner sep=2,minimum width=14,fill=white] (v90) at (0,-3.5) {$ b^* $};
	\draw[thick]  (v6) edge (v9);
	\draw [thick] (v90) edge (v2);
	\end{tikzpicture}
	\end{array}
	\end{align*}
	where $ b $ ranges over a basis of $ \Hom_{\C}(S,X^\vee T) $. We also have
	\begin{align*}
	d(S)
	\begin{array}{c}
	\begin{tikzpicture}[scale=0.15,every node/.style={inner sep=0,outer sep=-1}]
	\node (v5) at (-2,6) {};
	\node (v7) at (-2,-6) {};
	\node (v8) at (2,6) {};
	\node (v10) at (2,-6) {};
	\node at (-4,2) {$S$};
	\node at (3.5,2) {$X$};
	\draw[thick]  (v5) edge (v7);
	\draw[thick]  (v8) edge (v10);
	\end{tikzpicture}
	\end{array}
	= \sum\limits_{T,b} d(T)  \
	\begin{array}{c}
	\begin{tikzpicture}[scale=0.15,every node/.style={inner sep=0,outer sep=-1},yscale=1,xscale=-1]
	\node (v13) at (-0.5,1.5) {};
	\node (v130) at (-0.5,-2) {};
	\node (v6) at (0,7) {};
	\node (v2) at (0,-7.5) {};
	\node at (3,-6.5) {$S$};
	\node at (-5,-6.5) {$X$};
	\node at (3,6) {$S$};
	\node at (3,-0.25) {$T$};
	\node at (-5,6) {$X$};
	\draw [thick] plot[smooth, tension=.7] coordinates {(0.5441,1.6) (0.5441,-2.1)};
	\node (v1) at (-3,1.5) {};
	\node (v4) at (-3,-2) {};
	\node (v5) at (-3,-7.5) {};
	\node (v3) at (-3,7) {};
	\draw [thick] (v13) to[out=-90,in=-90] (v1);
	\draw [thick] (v3) edge (v1);
	\draw [thick] (v130) to[out=90,in=90] (v4);
	\draw [thick] (v5) edge (v4);
	\node (v7) at (-0.5,2.5) {};
	\node (v8) at (-0.5,-3) {};
	\draw [thick] (v7) edge (v13);
	\draw [thick] (v8) edge (v130);
	\node [draw,outer sep=0,inner sep=2,minimum width=14,fill=white] (v9) at (0,3) {$ b $};
	\node [draw,outer sep=0,inner sep=2,minimum width=14,fill=white] (v90) at (0,-3.5) {$ b^* $};
	\draw[thick]  (v6) edge (v9);
	\draw [thick] (v90) edge (v2);
	\end{tikzpicture}
	\end{array}
	\end{align*}
	where $ b $ ranges over a basis of $ \Hom_{\C}(S,T X^\vee) $.
	\proof The map $ \alpha(b) = \! $
	\begin{math}
	\begin{array}{c}
	\begin{tikzpicture}[scale=0.15,every node/.style={inner sep=0,outer sep=-1}]
	\node (v11) at (-3,1.5) {};
	\node (v13) at (-0.5,1.5) {};
	\node (v1) at (-3,7) {};
	\node (v6) at (0,7) {};
	\node at (-5,6) {$X$};
	\node at (2,0) {$T$};
	\node at (1.5,6) {$S$};
	\draw [thick] plot[smooth, tension=.7] coordinates {(0.5441,1.6) (0.5441,-1.1)};
	\draw [thick] (v13) to[out=-90,in=-90] (v11);
	\draw [thick] (v11) edge (v1);
	\node (v2) at (-0.5,2) {};
	\draw [thick] (v13) edge (v2);
	\node [draw,outer sep=0,inner sep=2,minimum size=10,fill=white] (v9) at (0,3) {$ b $};
	\draw[thick]  (v6) edge (v9);
	\end{tikzpicture}
	\end{array}
	\end{math}
	is the image of $ b $ under the adjunction
	\begin{align*}
	\alpha\colon \Hom_{\C}(S,X^\vee T) \to \Hom_{\C}(XS,T)
	\end{align*}
	and similarly the map $ \beta(b^*) \colon= \! $
	\begin{math}
	\begin{array}{c}
	\begin{tikzpicture}[scale=0.15,every node/.style={inner sep=0,outer sep=-1}]
	\node (v130) at (-0.5,-2) {};
	\node (v100) at (-3,-7.5) {};
	\node (v110) at (-3,-2) {};
	\node (v2) at (0,-7.5) {};
	\node at (-5,-6.5) {$X$};
	\node at (1.5,-6.5) {$S$};
	\node at (2,-0.5) {$T$};
	\draw [thick] plot[smooth, tension=.7] coordinates {(0.5441,0.8) (0.5441,-2.1)};
	\draw [thick] (v130) to[out=90,in=90] (v110);
	\draw [thick] (v110) edge (v100);
	\node (v1) at (-0.5,-3) {};
	\draw [thick] (v130) edge (v1);
	\node [draw,outer sep=0,inner sep=2,minimum size=10,fill=white] (v90) at (0,-3.5) {$ b^* $};
	\draw [thick] (v90) edge (v2);
	\end{tikzpicture}
	\end{array}
	\end{math}
	is the image of $ b^* $ under the adjunction
	\begin{align*}
	\beta\colon \Hom_{\C}(X^\vee T,S) \to \Hom_{\C}(T,XS).
	\end{align*}
	Therefore as $ b $ ranges over a basis of $ \Hom_{\C}(S,X^\vee T) $, $ \alpha(b) $ ranges over a basis $ \G $ of $ \Hom_{\C}(XS,T) $ and $ \beta(b^*) $ ranges over a basis $ \cH $ of $ \Hom_{\C}(T,XS) $. However $ \G $ and $ \cH $ are \emph{not} dual to one another. Indeed evaluating $ \alpha(b_1) \in \G $ on $ \beta(b_2^*) \in \cH $ gives
	\begin{align*}
	\begin{array}{c}
	\begin{tikzpicture}[scale=0.15,every node/.style={inner sep=0,outer sep=-1}]
	\node (v11) at (-2.75,1.2) {};
	\node (v13) at (-0.75,1.2) {};
	\node (v130) at (-0.75,10.75) {};
	\node (v110) at (-2.75,10.75) {};
	\node (v6) at (0,7) {};
	\node (v2) at (0,6.5) {};
	\node [draw,outer sep=0,inner sep=1,minimum height=14,minimum width=12] (v9) at (0,3) {$ b_1 $};
	\node [draw,outer sep=0,inner sep=1,minimum height=14,minimum width=12] (v90) at (0,9) {$ b_2^* $};
	\node at (-5,6) {$X$};
	\node at (2.8,0) {$T$};
	\node at (2.8,12) {$T$};
	\node at (2.8,6) {$S$};
	\draw[thick]  (v6) edge (v9);
	\draw [thick] (v90) edge (v2);
	\draw [thick] plot[smooth, tension=.7] coordinates {(0.5441,13) (0.5441,10.6)};
	\draw [thick] plot[smooth, tension=.7] coordinates {(0.5441,1.3) (0.5441,-1.1)};
	\draw [thick] (v110) edge (v11);
	\draw [thick] (v130) to[out=90,in=90] (v110);
	\draw [thick] (v13) to[out=-90,in=-90] (v11);
	\end{tikzpicture}
	\end{array}
	\refequals[Lem.~\ref{LEMMA:TWISTED_DUALS}] \frac{\tr(b_2^* \circ b_1)}{d(T)}
	\begin{array}{c}
	\begin{tikzpicture}[scale=0.15,every node/.style={inner sep=0,outer sep=-1}]
	\draw[thick] (0,7) -- (0,-7);
	\node at (2,4) {$ T $};
	\end{tikzpicture}
	\end{array} = \delta_{b_1,b_2} \frac{\tr(id_S)}{d(T)}
	\begin{array}{c}
	\begin{tikzpicture}[scale=0.15,every node/.style={inner sep=0,outer sep=-1}]
	\draw[thick] (0,7) -- (0,-7);
	\node at (2,4) {$ T $};
	\end{tikzpicture}
	\end{array} = \delta_{b_1,b_2} \frac{d(S)}{d(T)}
	\begin{array}{c}
	\begin{tikzpicture}[scale=0.15,every node/.style={inner sep=0,outer sep=-1}]
	\draw[thick] (0,7) -- (0,-7);
	\node at (2,4) {$ T $};
	\end{tikzpicture}
	\end{array}.
	\end{align*}
	This implies that $ \G $ and $ \cH $ are \emph{pseudo-dual}, to make them truly dual we would have to rescale one of them by $ \frac{d(T)}{d(S)} $. This, together with Lemma~\ref{LEMMA:DECOMPOSE}, proves the first part of the lemma, the second part is proved analogously.  \endproof
	\end{LEMMA}
An alternative proof of the above lemma may be found in \cite[Lemma 5.1]{MR2430629}.
	\begin{PROP}[Killing Ring] \label{PROP:KILLING_RING}
	Let $ R $ be in $ \Irr(\C) $. Then
		\begin{align*}
		\sum\limits_S d(S)
		\begin{array}{c}
		\begin{tikzpicture}[scale=0.15,every node/.style={inner sep=0,outer sep=-1}]
		\draw[thick] (-6,4) node (v4) {} to[out=90,in=90] (-3,4) node (v1) {};
		\draw[thick] (-3,-4) node (v3) {} to[out=-90,in=-90] (0,-4) node (v2) {};
		\draw[thick] (-2,1) -- (-3,0) -- (-6,-3) -- (-6,-6);
		\draw[thick] (-1,2) -- (0,3) -- (0,6);
		\draw [thick](v3) -- (-3,-3) -- (-4,-2) ;
		\draw [thick](v1) -- (-3,3) -- (0,0);
		\draw [thick](0,0) -- (0,-1.5) -- (v2);
		\node at (-7.5,2) {$ S $};
		\node at (2.2,-2) {$ S^\vee $};
		\node at (1.5,4.5) {$ R $};
		\draw[thick] (v4) -- (-6,0) -- (-5,-1);
		\end{tikzpicture}
		\end{array}
		= \delta_{R,\tid} d(\C)
		\end{align*}
	where $ \tid $ is the tensor identity.
	\proof See Corollary 3.1.11. in \cite{GraphCal}. \endproof
	\end{PROP}
	\begin{COR} \label{COR:KILLING_RING}
	Let $ R $ and $ S $ be in $ \Irr(\C) $. Then
		\begin{align*}
		\sum\limits_{S} d(S)
		\begin{array}{c}
		\begin{tikzpicture}[scale=0.15,every node/.style={inner sep=0,outer sep=-1}]
		\draw[thick] (-6,4) node (v4) {} to[out=90,in=90] (-3,4) node (v1) {};
		\draw[thick] (0,-7) node (v3) {} to[out=-90,in=-90] (3,-7) node (v2) {};
		\draw[thick] (1,-2) -- (0,-3) -- (-3,-6) -- (-3,-9);
		\draw[thick] (-2,1) -- (-3,0) -- (-6,-3) -- (-6,-9);
		\draw[thick] (2,-1) -- (3,0) -- (3,6);
		\draw[thick] (-1,2) -- (0,3) -- (0,6);
		\draw [thick] (-2,-4) -- (-4,-2) ;
		\draw [thick](v1) -- (-3,3) -- (3,-3);
		\draw [thick](3,-3) -- (3,-4.5) -- (v2);
		\node at (-7.5,2.5) {$ S $};
		\node at (5,-5) {$ S^\vee $};
		\node at (-1.5,5) {$ R $};
		\node at (4.5,5) {$ T $};
		\draw[thick] (v4) -- (-6,0) -- (-5,-1);
		\draw[thick] (-1,-5) -- (0,-6) -- (v3);
		\end{tikzpicture}
		\end{array}
		= \delta_{R^\vee,T} \frac{d(\C)}{d(R)}
		\begin{array}{c}
		\begin{tikzpicture}[scale=0.15,every node/.style={inner sep=0,outer sep=-1}]
		\draw[thick] (-6,-4.5) node (v4) {} to[out=90,in=90] (-3,-4.5) node (v1) {};
		\draw[thick] (-6,-1.5) node (v3) {} to[out=-90,in=-90] (-3,-1.5) node (v2) {};
		\node at (-1,0) {$ R^\vee $};
		\node at (-7.5,0) {$ R $};
		\node at (-7.5,-6) {$ R $};
		\node at (-1,-6) {$ R^\vee $};
		\draw[thick] (v4) -- (-6,-7.5);
		\draw [thick](v1) -- (-3,-7.5);
		\draw [thick](v3) -- (-6,1.5);
		\draw [thick](v2) -- (-3,1.5);
		\end{tikzpicture}
		\end{array}.
		\end{align*}
	\proof This follows immediately from Lemma~\ref{LEMMA:DECOMPOSE}, Proposition~\ref{PROP:KILLING_RING} and the fact that
		\begin{align*}
		\Hom_{\C}(RT,\tid) = \Hom_{\C}(T,R^\vee) = \delta_{R^\vee,T} \id_{R^\vee}.
		\end{align*}
	We note that the $ d(R)^{-1} $ term appears as the creation and annihilation morphisms are \emph{not} dual to one another, indeed they compose to the dimension. To make them dual we therefore weight by the inverse of the dimension. \endproof
	\end{COR}
	\begin{PROP} \label{PROP:TWISTED_S}
	We consider $ I,J,I',J' \in \Irr(\C) $, $ L $ in $ \C $, $ j \in \Hom_{\C}(IJ,L) $ and $ k \in \Hom_{\C}(L,I'J') $. Then
		\begin{align*}
		\sum\limits_{S} d(S)
		\begin{array}{c}
		\begin{tikzpicture}[scale=0.15,every node/.style={inner sep=0,outer sep=-1}]
		\node at (-7.5,2.5) {$ S $};
		\node at (5,-5) {$ S^\vee $};
		\node at (1.5,5) {$ I $};
		\node at (-8,-22) {$I' $};
		\node at (5.5,5) {$ J $};
		\node at (-3.5,-22) {$ J' $};
		\node at (0,-8.8) {$ L $};
		\node (v6) at (-2,-4.5) {};
		\node (v8) at (-1,-4.5) {};
		\node (v11) at (-2,-13) {};
		\node (v13) at (-1,-13) {};
		\draw[thick] (-6,4) node (v4) {} to[out=90,in=90] (-3,4) node (v1) {};
		\draw[thick] (0,-22) node (v3) {} to[out=-90,in=-90] (3,-22) node (v2) {};
		\draw [thick](v1) -- (-3,3) -- (1,-1);
		\draw [thick]  (2,-2) -- (3,-3) -- (3,-4.5) -- (v2);
		\draw[thick] (-1,-20) -- (0,-21) -- (v3);
		\draw [thick] (-2,-19) -- (-4,-17) ;
		\draw[thick] (v4) -- (-6,-15) -- (-4,-17);
		\draw[thick] (-2,1) -- (-3,0) -- (-3,-3) node (v5) {};
		\draw[thick] (4,6) -- (4,1) -- (2,-1) -- (1,-2) -- (0,-3) node (v7) {};
		\draw[thick] (-1,2) -- (0,3) -- (0,6);
		\draw[thick]  (v5) to[out=-90,in=90] (v6);
		\draw [thick] (v7) to[out=-135,in=90] (v8);
		\draw [thick](-5,-17) -- (-6,-18) -- (-6,-23);
		\draw [thick](-2,-23) -- (-2,-20) -- (0,-18) -- (0,-15) node (v14) {};
		\draw [thick] (v13) to[out=-90,in=90] (v14);
		\node (v12) at (-4,-16) {};
		\draw [thick] (v11) to[out=-90,in=45] (v12);
		\node [draw,outer sep=0,inner sep=1,minimum size=12, minimum height = 13, fill=white] (v9) at (-1.5,-5.5) {$ j $};
		\node [draw,outer sep=0,inner sep=1,minimum size=12, minimum height = 13, fill=white] (v10) at (-1.5,-12) {$ k $};
		\draw [thick] (v9) edge (v10);
		\end{tikzpicture}
		\end{array} = \delta_{I,I'} \delta_{J,J'} \tr(k \circ j)\frac{d(\C)}{d(I) d(J)}
		\begin{array}{c}
		\begin{tikzpicture}[scale=0.15,every node/.style={inner sep=0,outer sep=-1}]
		\node (v5) at (-2,14) {};
		\node (v7) at (-2,-14) {};
		\node (v8) at (2,14) {};
		\node (v10) at (2,-14) {};
		\node at (-1,12) {$I$};
		\node at (3.5,12) {$J$};
		\draw[thick]  (v5) edge (v7);
		\draw[thick]  (v8) edge (v10);
		\end{tikzpicture}
		\end{array}.
		\end{align*}
	\proof We have
		\begin{align*}
		&\sum\limits_{S} d(S)
			\begin{array}{c}
			\begin{tikzpicture}[scale=0.15,every node/.style={inner sep=0,outer sep=-1}]
			\node at (-7.5,2.5) {$ S $};
			\node at (5,-5) {$ S^\vee $};
			\node at (1.5,5) {$ I $};
			\node at (-8,-22) {$I' $};
			\node at (5.5,5) {$ J $};
			\node at (-3.5,-22) {$ J' $};
			\node at (0,-8.8) {$ L $};
			\node (v6) at (-2,-4.5) {};
			\node (v8) at (-1,-4.5) {};
			\node (v11) at (-2,-13) {};
			\node (v13) at (-1,-13) {};
			\draw[thick] (-6,4) node (v4) {} to[out=90,in=90] (-3,4) node (v1) {};
			\draw[thick] (0,-22) node (v3) {} to[out=-90,in=-90] (3,-22) node (v2) {};
			\draw [thick](v1) -- (-3,3) -- (1,-1);
			\draw [thick]  (2,-2) -- (3,-3) -- (3,-4.5) -- (v2);
			\draw[thick] (-1,-20) -- (0,-21) -- (v3);
			\draw [thick] (-2,-19) -- (-4,-17) ;
			\draw[thick] (v4) -- (-6,-15) -- (-4,-17);
			\draw[thick] (-2,1) -- (-3,0) -- (-3,-3) node (v5) {};
			\draw[thick] (4,6) -- (4,1) -- (2,-1) -- (1,-2) -- (0,-3) node (v7) {};
			\draw[thick] (-1,2) -- (0,3) -- (0,6);
			\draw[thick]  (v5) to[out=-90,in=90] (v6);
			\draw [thick] (v7) to[out=-135,in=90] (v8);
			\draw [thick](-5,-17) -- (-6,-18) -- (-6,-23);
			\draw [thick](-2,-23) -- (-2,-20) -- (0,-18) -- (0,-15) node (v14) {};
			\draw [thick] (v13) to[out=-90,in=90] (v14);
			\node (v12) at (-4,-16) {};
			\draw [thick] (v11) to[out=-90,in=45] (v12);
			\node [draw,outer sep=0,inner sep=1,minimum size=12, minimum height = 13, fill=white] (v9) at (-1.5,-5.5) {$ j $};
			\node [draw,outer sep=0,inner sep=1,minimum size=12, minimum height = 13, fill=white] (v10) at (-1.5,-12) {$ k $};
			\draw [thick] (v9) edge (v10);
			\end{tikzpicture}
			\end{array}
		= \sum\limits_{S} d(S)
			\begin{array}{c}
			\begin{tikzpicture}[scale=0.15,every node/.style={inner sep=0,outer sep=-1}]
			\node at (-7.5,2.5) {$ S $};
			\node at (8,-12) {$ S^\vee $};
			\node at (1.5,5) {$ I $};
			\node at (-7.5,-26) {$ I' $};
			\node at (5.5,5) {$ J $};
			\node at (-3.5,-26) {$ J' $};
			\node at (0,-8.8) {$ L $};
			\node [draw,outer sep=0,inner sep=1,minimum size=12, minimum height = 13 ] (v9) at (-1.5,-5.5) {$ j $};
			\node [draw,outer sep=0,inner sep=1,minimum size=12, minimum height = 13] (v10) at (-1.5,-12) {$ k $};
			\node (v6) at (-2,-3.8) {};
			\node (v8) at (-1,-3.8) {};
			\node (v11) at (-2,-13.68) {};
			\node (v13) at (-1,-13.68) {};
			\node (v80) at (2,-4) {};
			\draw[thick] (-6,4) node (v4) {} to[out=90,in=90] (-3,4) node (v1) {};
			\draw[thick] (2,-23) node (v3) {} to[out=-45,in=-90] (6,-23) node (v2) {};
			\draw [thick](v1) -- (-3,3) -- (1,-1);
			\draw [thick]  (1,-1) -- (6,-6)  -- (v2);
			\draw [thick] (-2,-19) -- (-4,-17) ;
			\draw[thick] (v4) -- (-6,-15) -- (-4,-17);
			\draw[thick] (-2,1) -- (-3,0) -- (-3,-3) node (v5) {};
			\draw[thick] (-1,2) -- (0,3) -- (0,6);
			\draw[thick]  (v5) to[out=-90,in=90] (v6);
			\draw [thick] (v80) to[out=90,in=90] (v8);
			\draw [thick] (v9) edge (v10);
			\draw [thick](-3,-15) node (v12) {} -- (-4,-16);
			\draw [thick](-5,-17) -- (-6,-18) -- (-6,-27);
			\draw [thick](-2,-27) -- (-2,-20) -- (0,-18) -- (0,-15) node (v14) {};
			\draw [thick] (v11) to[out=-90,in=-135] (v12);
			\draw [thick] (v13) to[out=-90,in=90] (v14);
			\draw [thick](v80) -- (2,-20) -- (0,-22) -- (0,-23) node (v7) {};
			\draw [thick](-1,-20) -- (0,-21);
			\draw [thick](1,-22) -- (2,-23);
			\draw [thick](4,6) -- (4,0) -- (10,-6) -- (10,-23) node (v15) {};
			\draw [thick] (v7) to[out = -90, in= -90] (v15);
			\end{tikzpicture}
			\end{array} \\
		= &\sum\limits_{S} d(S)
			\begin{array}{c}
			\begin{tikzpicture}[scale=0.15,every node/.style={inner sep=0,outer sep=-1}]
			\node at (5,-20) {$ S^\vee $};
			\node at (-0.5,-1) {$ I $};
			\node at (-7.5,-26) {$ I' $};
			\node at (5.5,-1) {$ J $};
			\node at (-3.5,-26) {$ J' $};
			\node at (0,-8.8) {$ L $};
			\node (v6) at (-2,-4.5) {};
			\node (v8) at (-1,-3.75) {};
			\node (v11) at (-2,-13.375) {};
			\node (v13) at (-1,-13.375) {};
			\node (v80) at (2,-3.75) {};
			\draw[thick]  (-2,0) to[out=-90,in=90] (v6);
			\draw [thick] (v80) to[out=90,in=90] (v8);
			\draw [thick](-3,-15) node (v12) {} -- (-5,-17);
			\draw [thick](-5,-17) -- (-6,-18) -- (-6,-27);
			\draw [thick](-2,-27) -- (-2,-20) -- (0,-18) -- (0,-15) node (v14) {};
			\draw [thick] (v11) to[out=-90,in=-135] (v12);
			\draw [thick] (v13) to[out=-90,in=90] (v14);
			\draw [thick](v80) -- (2,-20) -- (0,-22) -- (0,-23) node (v7) {};
			\draw [thick](-1,-20) -- (0,-21);
			\draw [thick](7,0) -- (7,-6) -- (7,-12) -- (7,-23) node (v15) {};
			\draw [thick] (v7) to[out = -90, in= -90] (v15);
			\draw[line width = 0.15cm,white] (-2,-19) node (v4) {} to[out=135,in=135] (0,-16) node (v1) {};
			\draw[thick] (-2,-19) node (v4) {} to[out=135,in=135] (0,-16) node (v1) {};
			\draw[thick] (1,-22) node (v3) {} to[out=-45,in=-45] (3,-19) node (v2) {};
			\draw [line width = 0.15cm,white]  (0,-16) to[out=-45,in=135] (3,-19);
			\draw [thick]  (0,-16) to[out=-45,in=135] (3,-19);
			\node (v5) at (-1,-4.5) {};
			\draw [thick] (v8) edge (v5);
			\node (v16) at (-2,-13) {};
			\node (v17) at (-1,-13) {};
			\draw [thick] (v11) edge (v16);
			\draw [thick] (v13) edge (v17);
			\node [draw,outer sep=0,inner sep=1,minimum size=12, minimum height = 13,fill=white] (v9) at (-1.5,-5.5) {$ j $};
			\node [draw,outer sep=0,inner sep=1,minimum size=12, minimum height = 13,fill=white] (v10) at (-1.5,-12) {$ k $};
			\draw [thick] (v9) edge (v10);
			\end{tikzpicture}
			\end{array}
		\refequals[Cor.~\ref{COR:KILLING_RING}] \delta_{J,J'} \frac{d(\C)}{d(J)}
			\begin{array}{c}
			\begin{tikzpicture}[scale=0.15,every node/.style={inner sep=0,outer sep=-1},yscale=1,xscale=-1]
			\node (v11) at (-3,1.25) {};
			\node (v13) at (-0.5,1.25) {};
			\node (v130) at (-0.5,13.75) {};
			\node (v110) at (-3,13.75) {};
			\node (v6) at (0,10) {};
			\node (v2) at (0,10) {};
			\node at (-5,7.5) {$J^\vee$};
			\node at (-8.5,7.5) {$J$};
			\node at (2,0) {$I'$};
			\node at (2.2,15) {$I$};
			\node at (1.5,7.5) {$L$};
			\draw [thick] plot[smooth, tension=.7] coordinates {(0.5,16) (0.5,13)};
			\draw [thick] plot[smooth, tension=.7] coordinates {(0.5,2) (0.5,-1)};
			\draw[thick] (-7,16) -- (-7,-1);
			\draw [thick] (v130) to[out=90,in=90] (v110);
			\draw [thick] (v13) to[out=-90,in=-90] (v11);
			\draw [thick] (v110) edge (v11);
			\node (v1) at (-0.5,13) {};
			\draw [thick] (v130) edge (v1);
			\node (v3) at (-0.5,2) {};
			\draw [thick] (v3) edge (v13);
			\node [draw,outer sep=0,inner sep=2.2,minimum size=13, minimum height = 13,fill=white] (v9) at (0,3) {$ k $};
			\node [draw,outer sep=0,inner sep=2.2,minimum size=13, minimum height = 13,fill=white] (v90) at (0,12) {$ j $};
			\draw[thick]  (v6) edge (v9);
			\draw [thick] (v90) edge (v2);
			\end{tikzpicture}
			\end{array} \\
		&\refequals[Lem.~\ref{LEMMA:TWISTED_DUALS}] \delta_{I,I'} \delta_{J,J'}\tr(k \circ j)\frac{d(\C)}{d(I) d(J)}
			\begin{array}{c}
			\begin{tikzpicture}[scale=0.15,every node/.style={inner sep=0,outer sep=-1}]
			\node (v5) at (-2,14) {};
			\node (v7) at (-2,0) {};
			\node (v8) at (2,14) {};
			\node (v10) at (2,0) {};
			\node at (-1,12) {$I$};
			\node at (3.5,12) {$J$};
			\draw[thick]  (v5) edge (v7);
			\draw[thick]  (v8) edge (v10);
			\end{tikzpicture}
			\end{array}.
			\end{align*}
	\endproof
	\end{PROP}

 \section{Introduction to $ \TC $} \label{SEC:INTROTC}

 We now introduce a new category, denoted $ \TC $, which shares the same objects as $ \C $ but admits more morphisms (i.e. $ \Hom_{\C}(X,Y) \leq \Hom_{\TC}(X,Y) $).  The intuition is that whereas morphisms in $ \C $ can be represented graphically as diagrams drawn on a bounded region of the plane, morphisms in $ \TC $ are given  by diagrams drawn on a \emph{cylinder}. For example, for any $ f \in \Hom_{\C}(X,Y) $ diagrammatically represented by
	\begin{align*}
	\begin{tikzpicture}[scale=0.15,every node/.style={inner sep=0,outer sep=-1}]
	\node (v2) at (4,6) {};
	\node (v1) at (-4,6) {};
	\node (v4) at (4,-6) {};
	\node (v3) at (-4,-6) {};
	\node (v5) at (0,6) {};
	\node (v7) at (0,-6) {};
	\node at (2,5) {$X$};
	\node at (2,-5) {$Y$};
	\node [draw,outer sep=0,inner sep=3,minimum size=10] (v6) at (0,0) {$ f $};
	\draw[thick]  (v6) edge (v5);
	\draw[thick]  (v7) edge (v6);
	\end{tikzpicture}
	\end{align*}
 there will be a morphism in $ \TC $ diagrammatically represented by
	\begin{align} \label{EQ:CYLINDER_HOM}
	\begin{array}{c}
	\begin{tikzpicture}[scale=0.15,every node/.style={inner sep=0,outer sep=-0.5}]
	\node (v2) at (4,12) {};
	\node (v1) at (-4,12) {};
	\node (v4) at (4,-6) {};
	\node (v3) at (-4,-6) {};
	\node (v5) at (0,10) {};
	\node (v7) at (-4,2) {};
	\node at (1,8.5) {$X$};
	\node at (0.5,-4.5) {$Y$};
	\node [draw,outer sep=0,inner sep=2,minimum size=10] (v6) at (-2,6) {$ f $};
	\draw[thick]  (v6) edge (v5);
	\draw[thick]  (v7) edge (v6);
	\draw  (0,12) ellipse (4 and 2);
	\draw  (v3) edge (v1);
	\draw  (v2) edge (v4);
	\draw (v3) to[out=-90,in=-90] (v4);
	\node [very thick] (v8) at (4,-1) {};
	\node [very thick] (v9) at (0,-8.2) {};
	\draw[thick,dotted]  (v7) to[in=90,out=-90] (v8);
	\draw[thick]  (v8) edge (v9);
	\end{tikzpicture}
	\end{array}.
	\end{align}
We capture such morphisms by drawing diagrams in a diamond and glueing the upper left and lower right edges. For example morphism \eqref{EQ:CYLINDER_HOM} is represented by
	\begin{align*}
	\begin{array}{c}
	\begin{tikzpicture}[scale=0.5,every node/.style={inner sep=0,outer sep=-1}]
	\node (v1) at (0,5) {};
	\node (v4) at (0,-5) {};
	\node (v2) at (5,0) {};
	\node (v3) at (-5,0) {};
	\node (v9) at (2.5,2.5) {};
	\node (v6) at (-2.5,-2.5) {};
	\node (v8) at (-3.5,1.5) {};
	\node (v11) at (1.5,-3.5) {};
	\node [draw,outer sep=0,rotate=-45,inner sep=2,minimum size=10] (v5) at (1,1) {\mbox{\Huge $ f $}};
	\node at (-3,-3) {$ Y $};
	\node at (3,3) {$ X $};
	\node at (2,-4) {$ Y^\vee $};
	\node at (-4,2) {$ Y^\vee $};
	\draw[thick]  (v9) edge (v5);
	\draw[thick]  plot[smooth, tension=.7] coordinates {(v8) (-1,0) (0.3031,0.3603)};
	\draw[thick]  plot[smooth, tension=.7] coordinates {(v11) (-0.5,-2) (v6)};
	\draw[very thick, red]  (v1) edge (v3);
	\draw[very thick, red]  (v2) edge (v4);
	\draw[very thick]  (v1) edge (v2);
	\draw[very thick]  (v3) edge (v4);
	\end{tikzpicture}
	\end{array}.
	\end{align*}
We note that this diagram can be read vertically as an element in $ \Hom_{\C}(Y^\vee X, Y Y^\vee) $. We also note that due to Lemma~\ref{LEMMA:DECOMPOSE} we may restrict ourselves to only gluing \emph{simple} stands. In this way morphism \eqref{EQ:CYLINDER_HOM} would be represented as
	\begin{align*}
	\sum\limits_{R,b}
	\begin{array}{c}
	\begin{tikzpicture}[scale=0.5,every node/.style={inner sep=0,outer sep=-1}]
	\node (v1) at (0,5) {};
	\node (v4) at (0,-5) {};
	\node (v2) at (5,0) {};
	\node (v3) at (-5,0) {};
	\node (v9) at (2.5,2.5) {};
	\node (v6) at (-2.5,-2.5) {};
	\node (v8) at (-2.208,0.5701) {};
	\node (v11) at (0.5,-2.375) {};
	\node (v7) at (-3.5,1.5) {};
	\node (v10) at (-2.9572,1.0566) {};
	\node at (-4,2) {$ R $};
	\node at (2,-3.75) {$ R $};
	\node [draw,outer sep=0,rotate=-45,inner sep=2,minimum size=10] (v5) at (1,1) {\mbox{\Huge $ f $}};
	\node at (-3,-3) {$ Y $};
	\node at (3,3) {$ X $};
	\node at (0.25,-1.5) {$ Y^\vee $};
	\node at (-1.125,0.75) {$ Y^\vee $};
	\draw[thick]  (v9) edge (v5);
	\draw[thick]  plot[smooth, tension=.7] coordinates {(v8) (-1,0) (0.3031,0.3603)};
	\draw[thick]  plot[smooth, tension=.7] coordinates {(v11) (-0.5,-2) (v6)};
	\draw[thick]  plot[smooth, tension=1] coordinates {(v7) (-3.1958,1.2318) (v10)};
	\draw[thick]  plot[smooth, tension=.7] coordinates {(0.5,-2.375) (1.3339,-3.1131) (1.6281,-3.3956)};
	\draw[very thick, red]  (v1) edge (v3);
	\draw[very thick, red]  (v2) edge (v4);
	\draw[very thick]  (v1) edge (v2);
	\draw[very thick]  (v3) edge (v4);
	\draw [thick] (v10) edge (v8);
	\node [draw,rotate=55,outer sep=0,inner sep=2,minimum size=14,fill=white] at (-2.5593,0.83) {$ b $};
	\node [draw,rotate=50,outer sep=0,inner sep=2,minimum size=14,fill=white] at (0.7416,-2.5955) {$ b^* $};
	\end{tikzpicture}
	\end{array}
	\end{align*}
where $ b $ ranges over a basis of $ \Hom_{\C}(R,Y^\vee) $. We note that each diagram may now be read vertically as an element in $ \Hom_{\C}(RX,YR) $. With this motivation in mind we may proceed with the definition of $ \TC $.
	\begin{DEF}
	Let $ \C $ be an MTC. The associated \emph{tube category}, denoted $ \TC $, is defined as the following category,
		\begin{enumerate}
		\item $ \Obj(\TC) \defeq \Obj(\C) $
		\item $ \Hom_{\TC}(X,Y) \defeq \bigoplus_{R} \Hom_{\C}(RX,YR) $
		\item Let $ f $ be in $ \Hom_{\TC}(X,Y)  $ and let $ g $ be in $ \Hom_{\TC}(Y,Z)  $. We define $ g \circ f $ as follows (using the diagrams explained above)
			\begin{align*}
			g \circ f \colon= \bigoplus\limits_{T} \sum\limits_{S,R,b}
			\begin{array}{c}
			\begin{tikzpicture}[scale=0.5,every node/.style={inner sep=0,outer sep=-1}]
			\node (v1) at (-1,4) {};
			\node (v4) at (-1.5,-6.5) {};
			\node (v2) at (4,-1) {};
			\node (v3) at (-6.5,-1.5) {};
			\node (v9) at (1.5,1.5) {};
			\node (v6) at (-4,-4) {};
			\node (v70) at (1.25,-3.75) {};
			\node (v7) at (-3.75,1.25) {};
			\node at (1.2,-1.6) {$ R $};
			\node at (-3.6,-0.8) {$ S $};
			\node at (2,2) {$ X $};
			\node at (-4.5,-4.5) {$ Z $};
			\node at (-1.5607,-1.0119) {$ Y $};
			\node at (1.75,-4.25) {$ T $};
			\node at (-4.25,1.75) {$ T $};
			\draw[thick]  plot[smooth, tension=.7] coordinates {(-3,-2.25) (-3.3997,-1.3889) (-2.625,-0.5) (-2.75,0)};
			\draw[thick]  plot[smooth, tension=.7] coordinates {(-2.25,-3) (-1.3889,-3.3997) (-0.586,-2.6724) (0,-2.75)};
			\draw[thick]  plot[smooth, tension=.7] coordinates {(0.5,-0.25) (0.925,-1.0642) (0.0727,-1.9921) (0.25,-2.5)};
			\draw[thick]  plot[smooth, tension=.7] coordinates {(-0.25,0.5) (-1.0642,0.925) (-1.9921,0.0727) (-2.5,0.25)};
			\node (v11) at (-3,-2.25) {};
			\node (v12) at (-2.25,-3) {};
			\draw [thick] (v11) edge (v12);
			\node (v14) at (-0.25,0.5) {};
			\node (v13) at (0.5,-0.25) {};
			\draw [thick] (v13) edge (v14);
			\node [diamond,draw,outer sep=0,inner sep=0.3,minimum size=10,fill=white] (v50) at (-2.5,-2.5) {\mbox{\Large $ g_{_S} $}};
			\node [diamond,draw,outer sep=0,inner sep=-0.2,minimum size=10,fill=white] (v5) at (0.125,0.125) {\mbox{\Large $ f_R $}};
			\node [draw,rotate=45,outer sep=0,inner sep=2,minimum size=10,fill=white] (v8) at (-2.8,0.3) {$ b $};
			\node [draw,rotate=45,outer sep=0,inner sep=2,minimum size=10,fill=white] (v10) at (0.3,-2.8) {$ b^* $};
			\draw[thick]  (v50) edge (v5);
			\draw [thick] (v50) edge (v6);
			\draw[thick]  (v7) edge (v8);
			\draw[thick]  (v10) edge (v70);
			\draw[thick]  (v9) edge (v5);
			\draw[very thick, red]  (v1) edge (v3);
			\draw[very thick, red]  (v2) edge (v4);
			\draw[very thick]  (v1) edge (v2);
			\draw[very thick]  (v3) edge (v4);
			\end{tikzpicture}
			\end{array}
			\\ \in \bigoplus\limits_{T} \Hom_{\C}(TX,ZT) = \Hom_{\TC}(X,Z)
			\end{align*}
		where $ f_R $ and $ g_S $ are the $ \Hom_{\C}(RX,YR) $ and $ \Hom_{\C}(SY,ZS) $ components of $ f $ and $ g $ respectively and $ b $ ranges over a basis of $ \Hom_{\C}(T,SR) $.
 		\end{enumerate}
	\end{DEF}

From Lemma~\ref{LEMMA:DECOMPOSE} we see that this definition agrees with the intuition that composition corresponds to vertically stacking the cylinders upon which the diagrams are drawn. This intuition, together with the associativity of the tensor product, makes it clear that composition in $ \TC $ is associative and that the identity in $ \End_{\TC}(X) = \bigoplus_{R}\Hom_{\C}(RX,XR) $ is given by the element with $ R $-component $ \delta_{\tid,R} \ \id_{\C}(X) $.

	\begin{REM}
	If we consider the algebra
		\begin{align*}
		\TA \defeq \End_{\TC}\left(\bigoplus\limits_{S} S \right)
		\end{align*}
	we recover Oceanu's \emph{tube algebra} \cite{Ocneanu1993}. As $ \bigoplus_{S} S $ is a projective generator in $ \TC $ the functor
		\begin{align*}
		\RTC &\to \RTA\\
		F &\mapsto \Hom_{\RTC}\left(F,\bigoplus\limits_{S} S\right)
		\end{align*}
	gives an equivalence, i.e. $ \TC $ is Morita equivalent to $ \TA $.
	\end{REM}
\section{$ \TC $ and $ \RTC $}\label{SECTCRTC}
We start by recalling that the centre of $ \C $, denoted $ Z(\C) $, is a category with objects $ (X,\tau) $ where $ X $ is in $ \C $ and $ \tau $ is a half braiding on $ X $ (see \cite[Section 7.13]{EtingofBook} for the precise definitions). Let $ \overline{\C} $ be the modular tensor category obtained by equipping $ \C $ with the opposite braiding. We consider the functor
	\begin{align*}
	\Psi_1 \colon \CC &\to Z(\C) \\
	X \boxtimes Y &\mapsto (X \otimes Y,\sigma_X \otimes \bar{\sigma}_Y)
	\end{align*}
where $ \boxtimes $ denotes the Deligne tensor product, and the functor
	\begin{align*}
	\Psi_2 \colon Z(\C) &\to \Rep(\TC) \\
	(X,\tau) &\mapsto F_{X,\tau}
	\end{align*}
	where  $ F_{X,\tau} $ is given by
		\begin{align*}
		F_{X,\tau}(Y) = \Hom_{\C}(Y,X)
		\end{align*}
	and, for $ f = \bigoplus_{S} \begin{array}{c}
	\begin{tikzpicture}[scale=0.15,every node/.style={inner sep=0,outer sep=-1}]
	\node at (3.5,3) {$Z$};
	\node at (-3,3) {$S$};
	\node at (-3.5,-3) {$Y$};
	\node at (3.5,-3) {$S$};
	\node [draw,diamond,outer sep=0,inner sep=2,minimum size=10] (v6) at (0,0) {$ f_S $};
	\node at (0.5,5) {\text{}};
	\end{tikzpicture}
	\end{array}
	\in \Hom_{\TC}(Z,Y) $,
		\begin{align*} F_{X,\tau}(f) \colon \Hom_{\C}(Y,X) &\to \Hom_{\C}(Z,X) \\
		\begin{array}{c}
		\begin{tikzpicture}[scale=0.15,every node/.style={inner sep=0,outer sep=-1}]
		\node (v2) at (4,6) {};
		\node (v1) at (-4,6) {};
		\node (v4) at (4,-6) {};
		\node (v3) at (-4,-6) {};
		\node (v5) at (0,6) {};
		\node (v7) at (0,-6) {};
		\node at (2,5) {$Y$};
		\node at (2,-5) {$X$};
		\node [draw,outer sep=0,inner sep=3,minimum size=10] (v6) at (0,0) {$ g $};
		\draw[thick]  (v6) edge (v5);
		\draw[thick]  (v7) edge (v6);
		\end{tikzpicture}
		\end{array} &\mapsto \sum\limits_{S}
		\begin{array}{c}
		\begin{tikzpicture}[scale=0.15,every node/.style={inner sep=0,outer sep=-1}]
		\node (v21) at (-6,-9.5) {};
		\node (v19) at (-2,-9.5) {};
		\node (v22) at (-6,-14.5) {};
		\node (v20) at (-2,-14.5) {};
		\node (v16) at (-2,-3) {};
		\node (v18) at (2,-3) {};
		\node (v17) at (2,3) {};
		\node (v15) at (-2,3) {};
		\draw [thick] (v15) edge (v16);
		\draw [thick] (v17) edge (v18);
		\draw [thick] (v19) edge (v20);
		\draw [thick] (v21) edge (v22);
		\node (v5) at (2,7) {};
		\node (v50) at (-2,4) {};
		\node (v7) at (2,-16) {};
		\node (v70) at (-2,-9.5) {};
		\node at (-6,4) (v8) {};
		\node at (3.5,6) {$Z$};
		\node at (-8,-2) {$S^\vee$};
		\node at (-7.5,-17) {$X$};
		\node at (-3.5,-8.7) {$X$};
		\node at (-3.5,-3) {$Y$};
		\node at (4,-8) {$S$};
		\node [draw,diamond,outer sep=0,inner sep=2,minimum size=10,fill=white] (v6) at (0,0) {$ f_S $};
		\node [draw,outer sep=0,inner sep=1,minimum size=10,fill=white] (v700) at (-2,-6) {$ g $};
		\draw[thick]  (v6)+(2,3) node (v11) {} edge (v5);
		\draw[thick]  (v6)+(-2,3) node (v3) {} edge (v50);
		\draw[thick]  (v6)+(-2,-3) node (v9) {} edge (v700);
		\draw[thick]  (v6)+(2,-3) node (v14) {} edge (v7);
		\draw[thick]  (v50) to[out=90,in=90] (v8);
		\node (v1) at (-6,-9.5) {};
		\draw [thick] (v8) edge (v1);
		\node (v4) at (-6,-14.5) {};
		\node[draw,outer sep=0,inner sep=2.3,minimum width=27,,fill=white] (v2) at (-4,-12) {$ \tau_{S^\vee} $};
		\node (v10) at (-6,-18) {};
		\draw [thick] (v4) edge (v10);
		\node (v12) at (-2,-14.5) {};
		\node (v13) at (-2,-16) {};
		\draw [thick] (v12) edge (v13);
		\draw [thick] (v13) to[out=-90,in=-90] (v7);
		\draw [thick] (v700) edge (v70);
		\end{tikzpicture}
		\end{array}.
	\end{align*}
	\begin{REM}
	A curious reader could check that $ F_{X,\tau} $ is a functor using Lemma~\ref{LEMMA:DECOMPOSE}.
	\end{REM}
\noindent We have the following result from \cite{EtingofBook} and \cite{Vaes}.
	\begin{THM} \label{THM:ETINGOF}
	$ \Psi_2 \circ \Psi_1 $ is a equivalence of categories from $ \CC $ to $ \RTC $.
	\proof Proposition 8.20.12. in \cite{EtingofBook} proves that $ \Psi_1 $ is an equivalence and Proposition 3.14 in \cite{Vaes} then proves that $ \Psi_2 $ is an equivalence.
 \endproof
	\end{THM}
For $ I, J \in \Irr(\C) $ we consider the functor $ F_{IJ} $, given by $ F_{IJ}(X) = \Hom_{\C}(X,IJ) $ for $ X $ in $ \C $ and
	\begin{align*}
	F_{IJ}(f) \colon \Hom_{\C}(Y,IJ) &\to \Hom_{\C}(Z,IJ)\\
	\begin{array}{c}
	\begin{tikzpicture}[scale=0.15,every node/.style={inner sep=0,outer sep=-1}]
	\node (v2) at (4,6) {};
	\node (v1) at (-4,6) {};
	\node (v4) at (4,-6) {};
	\node (v3) at (-4,-6) {};
	\node (v5) at (0,6) {};
	\node at (2,5) {$Y$};
	\node at (-2.25,-3.75) {$I$};
	\node at (2.25,-3.75) {$J$};
	\node [draw,outer sep=0,inner sep=3,minimum size=10] (v6) at (0,0) {$ g $};
	\draw[thick]  (v6) edge (v5);
	\node (v7) at (-0.75,-1.75) {};
	\node (v9) at (0.75,-1.75) {};
	\node (v8) at (-0.75,-6) {};
	\node (v10) at (0.75,-6) {};
	\draw [thick] (v7) edge (v8);
	\draw [thick] (v9) edge (v10);
	\end{tikzpicture}
	\end{array} &\mapsto \begin{array}{c}
	\begin{tikzpicture}[scale=0.15,every node/.style={inner sep=0,outer sep=-1}]
	\node (v30) at (-4,-5.5) {};
	\node (v50) at (-2,4) {};
	\node (v26) at (-2,0.5) {};
	\node (v26) at (-2,0.5) {};
	\node (v19) at (2,4) {};
	\node (v21) at (2,0.5) {};
	\node (v23) at (-2,-1.5) {};
	\node (v22) at (2,-4) {};
	\draw [thick] (v19) edge (v21);
	\draw [thick] (v22) edge (v21);
	\node (v24) at (-1.5,-1) {};
	\draw [thick] (v23) edge (v24);
	\draw [thick] (v50) edge (v26);
	\node (v5) at (2,7) {};
	\node (v7) at (2,-18) {};
	\node (v70) at (-6,-10) {};
	\node (v8) at (-10,4) {};
	\node (v14) at (-6,-6) {};
	\node at (3.5,6) {$Z$};
	\node at (-12,-2) {$S^\vee$};
	\node (v18) at (-7.5,-19) {$J$};
	\node (v18) at (-11.5,-19) {$I$};
	\node at (-4.5,-2) {$Y$};
	\node at (3.5,-8) {$S$};
	\node [draw,diamond,outer sep=0,inner sep=2,minimum size=10,fill=white] (v6) at (0,0) {$ f_S $};
	\draw[thick]  (v6)+(2,3) edge (v5);
	\draw[thick]  (v6)+(-2,3) node (v25) {} edge (v50);
	\draw[thick]  (v6)+(-2,-1.5) to[out=-135,in=90] (v30);
	\draw[thick]  (v6)+(2,-3) edge (v7);
	\draw[thick]  (v50) to[out=90,in=90] (v8);
	\node (v1) at (-10,-10) {};
	\draw [thick] (v8) edge (v1);
	\node (v3) at (-7.5,-11.5) {};
	\node (v9) at (-8.5,-12.5) {};
	\node (v4) at (-10,-14) {};
	\node (v2) at (-2,-6) {};
	\node (v20) at (-6,-14) {};
	\draw [thick] (v1) to[out=-90,in=90] (v20);
	\draw [thick] (v70) to[out=-90,in=45] (v3);
	\draw [thick] (v4) to[out=90,in=-135] (v9);
	\node (v10) at (-10,-20) {};
	\draw [thick] (v4) edge (v10);
	\node (v11) at (-6,-14) {};
	\node (v12) at (-2,-18) {};
	\draw [thick] (v11) to[out=-90,in=90] (v12);
	\node (v13) at (-2,-18) {};
	\draw [thick] (v13) to[out=-90,in=-90] (v7);
	\draw [thick] (v14) edge (v70);
	\node (v15) at (-2,-14) {};
	\draw [thick] (v2) edge (v15);
	\node (v16) at (-6,-18) {};
	\node (v17) at (-6,-20) {};
	\draw [line width= 0.15cm, white] (v15) to[out=-90,in=90] (v16);
	\draw [thick] (v15) to[out=-90,in=90] (v16);
	\draw [thick] (v16) edge (v17);
	\node [draw,outer sep=0,inner sep=1,minimum height=10,minimum width=25,fill=white] (v700) at (-4,-6) {$   g  $};
	\end{tikzpicture}
	\end{array}
	\end{align*}
for $ f = \bigoplus_{S} \begin{array}{c}
\begin{tikzpicture}[scale=0.15,every node/.style={inner sep=0,outer sep=-1}]
\node at (3.5,3) {$Z$};
\node at (-3,3) {$S$};
\node at (-3.5,-3) {$Y$};
\node at (3.5,-3) {$S$};
\node [draw,diamond,outer sep=0,inner sep=2,minimum size=10] (v6) at (0,0) {$ f_S $};
\node at (0.5,5) {\text{}};
\end{tikzpicture}
\end{array}
\in \Hom_{\TC}(Z,Y) $. \newpage
	\begin{COR} \label{COR:COREY_SIMPLE_FUNCTOR}
	The set
		\begin{align*}
		\{ F_{IJ} \}_{I,J \in \Irr(\C)}
		\end{align*}
	forms a complete set of simple objects in $ \RTC $.
	\proof This follows directly from the fact that $ \{I \boxtimes J \}_{I,J \in \Irr(\C)} $ forms a complete set of simple objects in $ \CC $ and $ \Psi_2 \circ \Psi_1(I \boxtimes J) = F_{IJ} $. \endproof
	\end{COR}
The following definition is inspired by the diagrammatic description of Ocneanu projections given in \cite[Section 2]{Evans1998}.
	\begin{DEF} \label{DEF:BILAM}
	For $ X, Y $ in $ \C $ and $ I,J \in \Irr(\C) $ we consider the (bi)linear map
		\begin{align*}
		\lambda_{YX}^{IJ}: \Hom_{\C}(IJ,Y) \otimes \Hom_{\C}(X,IJ) &\to \Hom_{\TC}(X,Y)\\
		\end{align*}
	given by
		\begin{align*}
		\lambda_{YX}^{IJ}(j \otimes i) \defeq \frac{d(I)d(J)}{d(\C)}\bigoplus\limits_{S} d(S)
		\begin{array}{c}
		\begin{tikzpicture}[scale=0.5,every node/.style={inner sep=0,outer sep=-1}]
		\node (v1) at (0,5) {};
		\node (v4) at (0,-5) {};
		\node (v2) at (5,0) {};
		\node (v3) at (-5,0) {};
		\node (v9) at (2.5,2.5) {};
		\node (v6) at (-2.5,-2.5) {};
		\node (v8) at (-2.5,2.5) {};
		\node (v11) at (2.5,-2.5) {};
		\node [every node] (v10) at (0.8,-0.8) {};
		\node [every node] (v12) at (1.2,-1.2) {};
		\node [every node] (v13) at (-1.5,-1.25) {};
		\node [every node] (v14) at (-1.2,0.8) {};
		\node [every node] (v15) at (-0.8,1.2) {};
		\node (v20) at (1.25,1.5) {};
		\node (v16) at (-1.25,-1.5) {};
		\node [every node] (v18) at (1,-1) {};
		\node (v17) at (1.5,1.25) {};
		\node [every node] at (3,3) {$ X $};
		\node [every node] at (-2.0068,-0.0314) {$ I $};
		\node [every node] at (1.9687,-0.0314) {$ J $};
		\node [every node] at (-3,-3) {$ Y $};
		\node [every node] at (-3,3) {$ S $};
		\node [every node] at (3,-3) {$ S $};
		\draw [thick] (v8) edge (v10);
		\draw [thick] (v11) edge (v12);
		\draw[thick]  plot[smooth, tension=.7] coordinates {(v20) (0.6,1) (-0.2,1.4) (v15)};
		\draw[thick]  plot[smooth, tension=.7] coordinates {(v13) (-1,-0.6) (-1.4,0.2) (v14)};
		\draw[thick]  plot[smooth, tension=.7] coordinates {(v16) (-0.6,-1) (0.2,-1.4) (v18) (1.4,-0.2) (1,0.6) (v17)};
		\node [draw,rotate=-45,outer sep=0,inner sep=1,minimum size=12,fill=white] (v5) at (1.5,1.5) {$i$};
		\node [draw,rotate=-45,outer sep=0,inner sep=1,minimum size=12,fill=white] (v7) at (-1.5,-1.5) {$j$};
		\draw[thick]  (v9) edge (v5);
		\draw [thick] (v6) edge (v7);
		\draw[very thick, red]  (v1) edge (v3);
		\draw[very thick, red]  (v2) edge (v4);
		\draw[very thick]  (v1) edge (v2);
		\draw[very thick]  (v3) edge (v4);
		\end{tikzpicture}
		\end{array}.
		\end{align*}
	\end{DEF}
	\begin{PROP}
	For $ j \in \Hom_{\C}(IJ,Y) $ the map
		\begin{align*}
		 \lambda_Y^{IJ}(j)_X: F_{IJ}(X) &\to \Hom_{\TC}(X,Y)\\
		i &\mapsto \lambda_{YX}^{IJ}(j \otimes i)
		\end{align*}
	is natural in $ X $, i.e. we have a map
		\begin{align*}
		\lambda_Y^{IJ} \colon \Hom_{\C}(IJ,Y) \to \Hom_\RTC(F_{IJ},Y^\sharp).
		\end{align*}
	for all $ Y $ in $ \C $.
	\proof We consider $ f = \bigoplus_{S} \begin{array}{c}
	\begin{tikzpicture}[scale=0.15,every node/.style={inner sep=0,outer sep=-1}]
	\node at (3.5,3) {$Z$};
	\node at (-3,3) {$S$};
	\node at (-3.5,-3) {$X$};
	\node at (3.5,-3) {$S$};
	\node [draw,diamond,outer sep=0,inner sep=2,minimum size=10] (v6) at (0,0) {$ f_S $};
	\node at (0.5,5) {\text{}};
	\end{tikzpicture}
	\end{array}
	\in \Hom_{\TC}(Z,X) $ and $ i \in \Hom_{\C}(X,IJ) $. We then compute \newpage
		\begin{align*}
		Y^\sharp(f) \circ (\lambda_Y^{IJ}(j))(i) &= \lambda_{YX}^{IJ}(j \otimes i)\circ f \\
		&= \frac{d(I)d(J)}{d(\C)}\bigoplus\limits_{R} \sum\limits_{S,T,b} d(T)
		\begin{array}{c}
		\begin{tikzpicture}[scale=0.5,every node/.style={inner sep=0,outer sep=-1}]
		\node (v1) at (1.5,6.5) {};
		\node (v4) at (1,-4) {};
		\node (v2) at (6.5,1.5) {};
		\node (v3) at (-4,1) {};
		\node (v6) at (-1.5,-1.5) {};
		\node (v8) at (-0.5,1.5) {};
		\node (v11) at (1.5,-0.5) {};
		\node (v29) at (4,4) {};
		\node (v28) at (-1,4) {};
		\node (v32) at (4,-1) {};
		\node [every node] (v10) at (0.75,0.25) {};
		\node [every node] (v12) at (1,0) {};
		\node [every node] (v13) at (-0.5,-0.25) {};
		\node [every node] (v14) at (0,1) {};
		\node [every node] (v15) at (0,1) {};
		\node (v20) at (1.25,1.5) {};
		\node (v16) at (-0.25,-0.5) {};
		\node [every node] (v18) at (0.75,-0.25) {};
		\node (v17) at (1.5,1.25) {};
		\node (v19) at (3,2) {};
		\node (v21) at (1.5,3.5) {};
		\node (v22) at (-0.375,1.375) {};
		\node (v23) at (0.25,0.25) {};
		\node (v24) at (1.5,-0.5) {};
		\node (v25) at (3.5,1.5) {};
		\node (v26) at (3,2) {};
		\node [every node] at (-0.75,0.5) {$ I $};
		\node [every node] at (0.625,-0.875) {$ J $};
		\node [every node] at (4.5,4.5) {$ Z $};
		\node [every node] at (2.25,-0.75) {$ T $};
		\node [every node] at (-1.5,4.5) {$ R $};
		\node [every node] at (4.5,-1.5) {$ R $};
		\node [every node] at (2,3.75) {$ S $};
		\node [every node] at (3.75,2) {$ S $};
		\node [every node] at (-2,-2) {$ Y $};
		\draw [thick] (v11) edge (v12);
		\draw [thick] (v26) edge (v25);
		\draw[thick]  plot[smooth, tension=.7] coordinates {(v20) (0.85,1.25) (0.25,1.25) (v15)};
		\draw[thick]  plot[smooth, tension=.7] coordinates {(v13) (-0.125,0.25) (-0.25,0.75) (v14)};
		\draw [line width=0.25cm,white] (v8)+(0.25,-0.25) edge (v10);
		\draw [thick] (v8) edge (v10);
		\draw[line width=0.25cm,white]  plot[smooth, tension=.7] coordinates { (0.25,-0.125) (v18) (1.25,0.25) (1.25,0.85) };
		\draw[thick]  plot[smooth, tension=.7] coordinates {(v16) (0.25,-0.125) (v18) (1.25,0.25) (1.25,0.85) (v17)};
		\draw [thick] (v19) edge (v21);
		\node (v27) at (0,2.75) {};
		\node (v31) at (0.25,3) {};
		\draw [thick]  (v8) to[out=135,in=-45] (v27);
		\draw [thick] (v21) to[out=135,in=-45] (v31);
		\node (v34) at (2.75,0) {};
		\node (v35) at (3,0.25) {};
		\draw [thick] (v11) to[out=-45,in=135] (v34);
		\draw [thick] (v25) to[out=-45,in=135] (v35);
		\node [draw,diamond,outer sep=0,inner sep=1,minimum size=25,fill=white] (v9) at (2.5,2.5) {$ f_S $};
		\node [draw,rotate=45,outer sep=0,inner sep=1,minimum size=12,fill=white] (v33) at (3.125,-0.125) {$ b^* $};
		\node [draw,rotate=45,outer sep=0,inner sep=1,minimum size=12,fill=white] (v30) at (-0.125,3.125) {$ b $};
		\node [draw,rotate=-45,outer sep=0,inner sep=1,minimum size=12,fill=white] (v5) at (1.5,1.5) {$ i $};
		\node [draw,rotate=-45,outer sep=0,inner sep=2,minimum size=12,fill=white] (v7) at (-0.5,-0.5) {$j$};
		\draw[thick]  (v9) edge (v5);
		\draw [thick] (v29) edge (v9);
		\draw [thick] (v6) edge (v7);
		\draw [thick] (v32) edge (v33);
		\draw [thick] (v28) edge (v30);
		\draw[very thick, red]  (v1) edge (v3);
		\draw[very thick, red]  (v2) edge (v4);
		\draw[very thick]  (v1) edge (v2);
		\draw[very thick]  (v3) edge (v4);
		\end{tikzpicture}
		\end{array}\\
		&= \frac{d(I)d(J)}{d(\C)}\bigoplus\limits_{R} \sum\limits_{S,T,b} d(T)
		\begin{array}{c}
		\begin{tikzpicture}[scale=0.5,every node/.style={inner sep=0,outer sep=-1}]
		\node (v1) at (1.5,6.5) {};
		\node (v4) at (0,-5) {};
		\node (v2) at (6.5,1.5) {};
		\node (v3) at (-5,0) {};
		\node (v6) at (-2.5,-2.5) {};
		\node (v29) at (4,4) {};
		\node (v28) at (-2.5,2.5) {};
		\node (v32) at (2.5,-2.5) {};
		\node [every node] (v10) at (-2,2) {};
		\node (v330) at (0.625,-0.625) {};
		\node [every node] (v13) at (-1.5,-1.25) {};
		\node [every node] (v14) at (-1.75,1.375) {};
		\node [every node] (v15) at (-1.375,1.75) {};
		\node (v20) at (1.25,1.5) {};
		\node (v16) at (-1.25,-1.5) {};
		\node [every node] (v18) at (0.25,-2) {};
		\node (v17) at (1.5,1.25) {};
		\node (v19) at (3,2) {};
		\node (v21) at (0.125,-0.5) {};
		\node (v22) at (-1,0.625) {};
		\node (v23) at (-0.625,1) {};
		\node (v26) at (3,2) {};
		\node (v31) at (0.5,-0.25) {};
		\node [every node] at (-2.375,0.625) {$ I $};
		\node [every node] at (1,-2.5) {$ J $};
		\node [every node] at (4.5,4.5) {$ Z $};
		\node [every node] at (-0.75,-0.375) {$ T $};
		\node [every node] at (3,-3) {$ R $};
		\node [every node] at (-3,3) {$ R $};
		\node [every node] at (3.75,2) {$ S $};
		\node [every node] at (2,3.75) {$ S $};
		\node [every node] at (-3,-3) {$ Y $};
		\draw[thick]  plot[smooth, tension=.7] coordinates {(v20) (0.85,1.25) (-0.25,2) (v15)};
		\draw[thick]  plot[smooth, tension=.7] coordinates {(v13) (-1.25,-0.625) (-2,0.25) (v14)};
		\draw [thick] (v32) edge (v330);
		\draw [thick] plot[smooth, tension=.7] coordinates {(v31) (0.125,0.375) (0.75,0.5) (2.375,-0.875) (3.125,-0.875) (4,0) (4,0.875) (v26)};
		\draw[line width=0.2cm,white]  plot[smooth, tension=.7] coordinates { (-0.625,-1.25) (v18) (1.5,-1.5) (2,-0.25) (1.25,0.85) };
		\draw[thick]  plot[smooth, tension=.7] coordinates {(v16) (-0.625,-1.25) (v18) (1.5,-1.5) (2,-0.25) (1.25,0.85) (v17)};
		\draw [thick] (v22) edge (v21);
		\draw [line width=0.2cm, white] plot[smooth, tension=.7] coordinates {(v23) (-0.125,0.625) (0.125,1.125) (-0.875,2.375) (-0.875,3.125) (-0.25,4) (0.875,4)};
		\draw [thick] plot[smooth, tension=.7] coordinates {(v23) (-0.125,0.625) (0.125,1.125) (-0.875,2.375) (-0.875,3.125) (-0.25,4) (0.875,4) (v19)};
		\node [draw,diamond,outer sep=0,inner sep=1,minimum size=25,fill=white] (v9) at (2.5,2.5) {$ f_S $};
		\node [draw,rotate=45,outer sep=0,inner sep=1,minimum size=12,fill=white] (v33) at (0.5,-0.5) {$ b^* $};
		\node [draw,rotate=45,outer sep=0,inner sep=1,minimum size=12,fill=white] (v30) at (-1,1) {$ b $};
		\node [draw,rotate=-45,outer sep=0,inner sep=1,minimum size=12,fill=white] (v5) at (1.5,1.5) {$ i $};
		\node [draw,rotate=-45,outer sep=0,inner sep=2,minimum size=12,fill=white] (v7) at (-1.5,-1.5) {$j$};
		\draw[thick]  (v9) edge (v5);
		\draw [thick] (v6) edge (v7);
		\draw [thick] (v29) edge (v9);
		\draw [thick] (v28) edge (v30);
		\draw[very thick, red]  (v1) edge (v3);
		\draw[very thick, red]  (v2) edge (v4);
		\draw[very thick]  (v1) edge (v2);
		\draw[very thick]  (v3) edge (v4);
		\end{tikzpicture}
		\end{array} \\
		&\refequals[Lem.~\ref{LEMMA:DUAL_DECOMPOSE}] \frac{d(I)d(J)}{d(\C)} \bigoplus\limits_{R} d(R) \sum\limits_{S}
		\begin{array}{c}
		\begin{tikzpicture}[scale=0.5,every node/.style={inner sep=0,outer sep=-1}]
		\node (v1) at (1.5,6.5) {};
		\node (v4) at (0,-5) {};
		\node (v2) at (6.5,1.5) {};
		\node (v3) at (-5,0) {};
		\node (v6) at (-2.5,-2.5) {};
		\node (v8) at (-3,2) {};
		\node (v11) at (2,-3) {};
		\node (v29) at (4,4) {};
		\node [every node] (v10) at (-0.25,-0.75) {};
		\node [every node] (v12) at (0,-1) {};
		\node [every node] (v13) at (-1.5,-1.25) {};
		\node [every node] (v14) at (-0.5,0.5) {};
		\node [every node] (v15) at (-0.5,0.5) {};
		\node (v20) at (1.25,1.5) {};
		\node (v16) at (-1.25,-1.5) {};
		\node [every node] (v18) at (0.45,-0.45) {};
		\node (v17) at (1.5,1.25) {};
		\node (v19) at (1.75,3) {};
		\node (v21) at (0,2.75) {};
		\node (v22) at (-0.375,1.375) {};
		\node (v23) at (0.25,0.25) {};
		\node (v24) at (1.375,-0.375) {};
		\node (v25) at (2.75,0) {};
		\node (v26) at (3,1.75) {};
		\node [every node] at (-1.5,-0.5) {$ I $};
		\node [every node] at (-0.5,-1.5) {$ J $};
		\node [every node] at (4.5,4.5) {$ Z $};
		\node [every node] at (2.5,-3.5) {$ R $};
		\node [every node] at (2.5,-3.5) {$ R $};
		\node [every node] at (3.625,0.875) {$ S $};
		\node [every node] at (-3,-3) {$ Y $};
		\draw [thick] (v11) edge (v12);
		\draw[thick]  plot[smooth, tension=.7] coordinates {(v20) (0.85,1.25) (0.25,1.25) (v15)};
		\draw[thick]  plot[smooth, tension=.7] coordinates {(v13) (-1,-0.6) (-0.75,0) (v14)};
		\draw [line width=0.25cm,white] (v8)+(1,-1) edge (v10);
		\draw [thick] (v8) edge (v10);
		\draw[line width=0.25cm,white]  plot[smooth, tension=.7] coordinates {(v19) (v21) (v22) (v23) (v24) (v25) (v26) };
		\draw[thick]  plot[smooth, tension=.7] coordinates {(v19) (v21) (v22) (v23) (v24) (v25) (v26)};
		\draw[line width=0.25cm,white]  plot[smooth, tension=.7] coordinates { (-0.6,-1) (0,-0.75) (v18) (1.25,0.25) (1.25,0.85) };
		\draw[thick]  plot[smooth, tension=.7] coordinates {(v16) (-0.6,-1) (0,-0.75) (v18) (1.25,0.25) (1.25,0.85) (v17)};
		\node (v27) at (2.5,2.5) {};
		\node (v28) at (2.5,2.5) {};
		\draw [thick] (v19) to[in=135,out=0] (v27);
		\draw [thick] (v26) to[in=-45,out=90] (v28);
		\node [draw,diamond,outer sep=0,inner sep=1,minimum size=25,fill=white] (v9) at (2.5,2.5) {$ f_S $};
		\node [draw,rotate=-45,outer sep=0,inner sep=1,minimum size=12,fill=white] (v5) at (1.5,1.5) {$ i $};
		\node [draw,rotate=-45,outer sep=0,inner sep=2,minimum size=12,fill=white] (v7) at (-1.5,-1.5) {$j$};
		\draw[thick]  (v9) edge (v5);
		\draw [thick] (v6) edge (v7);
		\draw [thick] (v29) edge (v9);
		\draw[very thick, red]  (v1) edge (v3);
		\draw[very thick, red]  (v2) edge (v4);
		\draw[very thick]  (v1) edge (v2);
		\draw[very thick]  (v3) edge (v4);
		\end{tikzpicture}
		\end{array}\\
		&=\lambda_{ZX}^{IJ}(j \otimes F_{IJ}(f)(i))  = \lambda_Z^{IJ}(j) \circ (F_{IJ}(f))(i),
		\end{align*}
	as required. \endproof
	\end{PROP}
\newpage We now consider $ \mu_Y^{IJ} : F_{IJ}(Y) \lra{\cong} \Hom_{\RTC}(Y^\sharp,F_{IJ}) $ as defined in Section~\ref{SEC:PRELIMSONYON} with $ F = F_{IJ} $.
	\begin{THM}
	 Under the trace pairing (c.f. Proposition~\ref{PROP:NONSIMPLE_PERFECT}), $ \lambda_Y^{IJ} $ is an opposite of $ \mu_Y^{IJ} $ (c.f. Definition~\ref{DEF:OPPOSITE}).
	\proof We consider $ X $ in $ \C $, $ j \in \Hom_{\C}(IJ,Y) $, $ k \in \Hom_{\C}(Y,IJ) $ and compute
		\begin{align*}
		(\mu_Y(k) \circ \lambda_Y(j))_X(i) &= \mu_Y(k)_X\left( \lambda_{YX}^{IJ}(j \otimes i) \right) = F_{IJ}(\lambda_{YX}^{IJ}(j \otimes i))(k)\\
		&= \frac{d(I)d(J)}{d(\C)}\bigoplus\limits_{S} d(S)
		\begin{array}{c}
		\begin{tikzpicture}[scale=0.5,every node/.style={inner sep=0,outer sep=-1}]
		\node (v9) at (2.5,3) {};
		\node (v8) at (-1.5,1.5) {};
		\node(v55) at (1.625,1.625) {};
		\node [every node] (v10) at (0.8,-0.8) {};
		\node [every node] (v12) at (1.2,-1.2) {};
		\node [every node] (v13) at (-1.5,-1.25) {};
		\node [every node] (v14) at (-1.2,0.8) {};
		\node [every node] (v15) at (-0.8,1.2) {};
		\node (v20) at (1.25,1.5) {};
		\node (v16) at (-1.25,-1.5) {};
		\node [every node] (v18) at (1,-1) {};
		\node (v17) at (1.5,1.25) {};
		\node [every node] at (1.6941,2.5345) {$ X $};
		\node [every node] at (-2.0068,-0.0314) {$ I $};
		\node [every node] at (-4,-5.75) {$ I $};
		\node [every node] at (-1,-5) {$ J $};
		\node [every node] at (1.9687,-0.0314) {$ J $};
		\node [every node] at (-2.5,-2.25) {$ Y $};
		\node [every node] at (2.25,-2.5) {$ S $};
		\draw[thick]  (v9) to[out=-90,in=45] (v55);
		\draw [thick] (v8) edge (v10);
		\draw[thick]  plot[smooth, tension=.7] coordinates {(v20) (0.6,1) (-0.2,1.4) (v15)};
		\draw[thick]  plot[smooth, tension=.7] coordinates {(v13) (-1,-0.6) (-1.4,0.2) (v14)};
		\draw[thick]  plot[smooth, tension=.7] coordinates {(v16) (-0.6,-1) (0.2,-1.4) (v18) (1.4,-0.2) (1,0.6) (v17)};
		\node (v3) at (-3,-4.5) {};
		\node (v4) at (-2.125,-3.625) {};
		\node (v19) at (-1.875,-3.625) {};
		\node (v6) at (-2.5,-4.5) {};
		\node (v21) at (-1.5,-4.5) {};
		\draw[thick]  (v8) to[out=135,in=135] (v3);
		\draw [thick] (v4) to[out=-90,in=45] (v6);
		\draw [thick] (v19) to[out=-90,in=90] (v21);
		\node (v22) at (-1.875,-5.625) {};
		\draw [thick] (v3) edge (v22);
		\node (v25) at (-3,-5) {};
		\node (v26) at (-3.5,-6.5) {};
		\node (v27) at (-1.375,-6.125) {};
		\node (v23) at (-1.5,-5.5) {};
		\node (v24) at (-2,-6.5) {};
		\draw [thick] (v21) edge (v23);
		\draw [thick] (v23) to[out=-90,in=90] (v24);
		\draw [thick] (v25) to[out=-135,in=90] (v26);
		\draw [thick] (v27) to[out=-45,in=-45] (v12);
		\node [draw,rotate=-45,outer sep=0,inner sep=1,minimum size=12,fill=white] (v5) at (1.5,1.5) {$i$};
		\node [draw,rotate=-45,outer sep=0,inner sep=1,minimum size=12,fill=white] (v7) at (-1.5,-1.5) {$j$};
		\node [draw,outer sep=0,inner sep=1,minimum size=12,fill=white] (v1) at (-2,-3.5) {$ k $};
		\draw[thick]  (v7) to[out=-110, in=90] (v1);
		\end{tikzpicture}
		\end{array}\\
		&= \tr(j \circ k) i
		\end{align*}
	where the last equality is due to Proposition~\ref{PROP:TWISTED_S}.
      Therefore $ \lambda_Y^{IJ} $ is an opposite of $ \mu_Y^{IJ} $. \endproof
	\end{THM}

Note that, by Corollary~\ref{COR:COREY_SIMPLE_FUNCTOR}, $ F_{IJ} $ is simple (and Schurian, as $ \fld $ is algebraically closed)
so  Remark~\ref{REM:SCHUR} applies, i.e. $ \lambda_Y^{IJ} $ is the dual of $ \bigl( \mu_Y^{IJ} \bigr)^{-1}$.
It now follows from Lemma~\ref{lem:compo}, that the (bi)linear version $ \lambda_{YX}^{IJ} $,
as in Definition~\ref{DEF:BILAM}, is also the map of the same name described in the introduction. \\

We can use this result to identify the \emph{primitive idempotents} in $ \End_{\TC}(X) $,
i.e. those that cannot be written as a non trivial direct sum.
	\begin{COR} \label{COR:PRIMITIVE}
	Let $ i \in \Hom_{\C}(X,IJ) $ and $ j \in \Hom_{\C}(IJ,X) $ be such that $ \tr(i \circ j) = 1 $. Then $ \idem = \lambda_{XX}^{IJ}(j \otimes i) $ is a primitive idempotent in $ \End_{\TC}(X) $.
	\proof By Corollary~\ref{COR:isoidem} we have $ (X,\idem)^\sharp \cong F_{IJ} $. By Corollary~\ref{COR:COREY_SIMPLE_FUNCTOR} $ F_{IJ} $ is simple, and therefore, $ \idem $ is primitive. \endproof
	\end{COR}

Finally, we can also enhance the composition formula in Proposition~\ref{PROP:compo} to reflect the fact that different $F_{IJ}$ are orthogonal.

	\begin{PROP} \label{PROP:compo2}
	We have the following composition rule,
	\begin{align*}
	\lambda_{ZY}^{I'J'}(l \otimes k) \circ \lambda_{YX}^{IJ}(j \otimes i) = \delta_{I,I'} \delta_{J,J'} \tr(k \circ j) \lambda_{ZX}^{IJ}(l \otimes i).
	\end{align*}
	\proof The case when $ I = I' $ and $ J= J' $ is given by Proposition~\ref{PROP:compo}. In the other case we consider $ \idem = \lambda_{XX}^{IJ}(i\dual \otimes i) $ and $ \idem'= \lambda_{ZZ}^{I'J'}(k\dual \otimes k) $ where $ i\dual \in \Hom_{\C}(IJ,X) $ and $ j\dual \in \Hom_{\C}(Z,IJ) $ are such that $ \tr(i \circ i\dual) = \tr(j \circ j\dual) = 1 $. As, once again by Proposition~\ref{PROP:compo},
		\begin{align*}
		\lambda_{ZY}^{I'J'}(l \otimes k) \circ \lambda_{YX}^{IJ}(j \otimes i) \circ \idem = \lambda_{ZY}^{I'J'}(l \otimes k) \circ \lambda_{YX}^{IJ}(j \otimes i)
		\end{align*}
	and
		\begin{align*}
		\idem' \circ 	\lambda_{ZY}^{I'J'}(l \otimes k) \circ \lambda_{YX}^{IJ}(j \otimes i) = 	\lambda_{ZY}^{I'J'}(l \otimes k) \circ \lambda_{YX}^{IJ}(j \otimes i)
		\end{align*}
	we have
		\begin{align*}
		\lambda_{ZY}^{I'J'}(l \otimes k) \circ \lambda_{YX}^{IJ}(j \otimes i) \circ \Blank \in \Hom_{\RTC}((X,\idem)^\sharp,(Z,\idem')^\sharp) = \Hom_{\RTC}(F_{IJ},F_{I'J'}) = 0,
		\end{align*}
	 where the last equality is due to Corollary~\ref{COR:COREY_SIMPLE_FUNCTOR}.
	\end{PROP}

	\begin{REM}
	The following graphical computation provides an alternative proof of Proposition~\ref{PROP:compo2} that does not use Corollary~\ref{COR:COREY_SIMPLE_FUNCTOR}. \\\\
		\begin{align*}
		&\lambda_{ZY}^{I'J'}(l \otimes k) \circ \lambda_{YX}^{IJ}(j \otimes i)  \\
		&=\frac{d(I)d(J)d(I')d(J')}{d(\C)^2} \bigoplus\limits_{R} \sum\limits_{S,T,b} d(S)d(T)
		\begin{array}{c}
		\begin{tikzpicture}[scale=0.5,every node/.style={inner sep=0,outer sep=-1}]
		\node at (0,0) {};
		\node (v1) at (0.25,5.25) {};
		\node (v4) at (-0.25,-5.25) {};
		\node (v2) at (5.25,0.25) {};
		\node (v3) at (-5.25,-0.25) {};
		\node (v9) at (2.75,2.75) {};
		\node (v6) at (-2.75,-2.75) {};
		\node (v70) at (2.5,-2.5) {};
		\node (v7) at (-2.5,2.5) {};
		\node at (2.25,0) {$ T $};
		\node at (-2.35,0.45) {$ S $};
		\node at (3.25,3.25) {$ X $};
		\node at (-3.25,-3.25) {$ Z $};
		\node at (3,-3) {$ R $};
		\node at (-3,3) {$ R $};
		\draw[thick]  plot[smooth, tension=.7] coordinates {(1.625,0.625) (1.7502,0.146) (1.3227,-0.7421) (1.5,-1.375)};
		\draw[thick]  plot[smooth, tension=.7] coordinates {(-0.625,-1.625) (-0.146,-1.7502) (0.7421,-1.3227) (1.375,-1.5)};
		\draw[thick]  plot[smooth, tension=.7] coordinates {(1.375,0.875) (0.375,1.75) (-0.25,1.75) (-0.75,1.25) (-1.25,1.5)};
		\draw[thick]  plot[smooth, tension=.7] coordinates {(-0.875,-1.375) (-1.75,-0.375) (-1.75,0.25) (-1.25,0.75) (-1.5,1.25)};
		\draw [thick] plot[smooth, tension=.7] coordinates {(1.625,1.5) (1.5,0.75) (0.75,0.5)};
		\draw [thick] plot[smooth, tension=.7] coordinates {(1.5,1.625) (1.1139,1.6049) (0.8822,1.5522)};
		\draw [thick] plot[smooth, tension=.7] coordinates {(0.5,0.75) (0.6,1.1) (0.6721,1.3221)};
		\draw [thick] plot[smooth, tension=.7] coordinates {(-1.5,-1.625) (-0.75,-1.5) (-0.625,-0.75)};
		\draw [thick] plot[smooth, tension=.7] coordinates {(-1.625,-1.5) (-1.625,-1.125) (-1.5522,-0.8822)};
		\draw [thick] plot[smooth, tension=.7] coordinates {(-0.75,-0.625) (-1.1,-0.6) (-1.3221,-0.6721)};
		\node [draw,outer sep=0,inner sep=1,rotate=-45,minimum size=12,fill=white] (v12) at (1.75,1.75) {$i$};
		\node [draw,outer sep=0,inner sep=1,rotate=-45,minimum size=12,fill=white] (v11) at (0.5,0.5) {$j$};
		\node [draw,outer sep=0,inner sep=2,rotate=-45,minimum size=12,fill=white] (v5) at (-0.5,-0.5) {$k$};
		\node [draw,outer sep=0,inner sep=2,rotate=-45,minimum size=12,fill=white] (v13) at (-1.75,-1.75) {$l$};
		\draw[thick]  (v9) edge (v12);
		\draw [thick] (v5) edge (v11);
		\draw [thick] (v13) edge (v6);
		\node [draw,rotate=45,outer sep=0,inner sep=2,minimum size=10,fill=white] (v8) at (-1.55,1.55) {$ b $};
		\node [draw,rotate=45,outer sep=0,inner sep=2,minimum size=10,fill=white] (v10) at (1.55,-1.55) {$ b^* $};
		\draw[thick]  (v7) edge (v8);
		\draw[thick]  (v10) edge (v70);
		\draw[very thick, red]  (v1) edge (v3);
		\draw[very thick, red]  (v2) edge (v4);
		\draw[very thick]  (v1) edge (v2);
		\draw[very thick]  (v3) edge (v4);
		\end{tikzpicture}
		\end{array} \\
		&= \frac{d(I)d(J)d(I')d(J')}{d(\C)^2} \bigoplus\limits_{R} \sum\limits_{S,T,b} d(S)d(T)
		\begin{array}{c}
		\begin{tikzpicture}[scale=0.5,every node/.style={inner sep=0,outer sep=-1}]
		\node (v1) at (2.25,7.25) {};
		\node (v4) at (-0.25,-5.25) {};
		\node (v2) at (7.25,2.25) {};
		\node (v3) at (-5.25,-0.25) {};
		\node (v9) at (4.75,4.75) {};
		\node (v6) at (-2.75,-2.75) {};
		\node[outer sep=-0.5] (v7) at (4.625,-0.375) {};
		\node[outer sep=-0.5] (v14) at (-0.375,4.625) {};
		\node at (4,2.5) {$ J $};
		\node at (2.5,4) {$ I $};
		\node at (2.5,2.625) {$ T $};
		\node at (-1,1.5) {$ S $};
		\node at (1.5,-1) {$ S^\vee $};
		\node at (5.25,5.25) {$ X $};
		\node at (-3.25,-3.25) {$ Z $};
		\node at (5.125,-0.875) {$ R $};
		\node at (-0.875,5.125) {$ R $};
		\node (v100) at (2.75,1.5) {};
		\draw [thick] (v100) edge (v7);
		\draw [thick] plot[smooth, tension=.7] coordinates {(3.4,3.625) (2.75,3.25) (1.5,3.625) (0.625,2.875) (1,1.5) (0.5,0.75)};
		\draw[thick]  plot[smooth, tension=.5] coordinates {(2.25,1.625) (1.9193,1.8143) (1.62,1.49) (2.125,0.375) (0,-1.625) (-0.625,-1.625)};
		\draw[line width=0.18cm,white]  plot[smooth, tension=.5] coordinates { (1.4128,1.7578) (0.375,2.125) (-1.625,0) (-1.625,-0.625)};
		\draw[thick]  plot[smooth, tension=.5] coordinates {(1.5,2.375) (1.7371,2.0571) (1.4128,1.7578) (0.375,2.125) (-1.625,0) (-1.625,-0.625)};
		\draw [thick] plot[smooth, tension=.7] coordinates {(-1.5,-1.625) (-0.75,-1.5) (-0.625,-0.75)};
		\draw [thick] plot[smooth, tension=.7] coordinates {(-1.625,-1.5) (-1.625,-1.125) (-1.5522,-0.8822)};
		\draw [thick] plot[smooth, tension=.7] coordinates {(-0.75,-0.625) (-1.1,-0.6) (-1.3221,-0.6721)};
		\draw [thick] plot[smooth, tension=.7] coordinates {(-1.625,-0.625) (-0.875,-1.375)};
		\draw [thick] plot[smooth, tension=.7] coordinates {(1.75,2.625)};
		\draw [thick] plot[smooth, tension=.7] coordinates {(1.75,2.625) (2.5,1.875)};
		\draw [line width=0.18cm,white] plot[smooth, tension=.7] coordinates { (3.25,2.75) (3.625,1.5) (2.875,0.625) (2.015,0.86) (1.5,1)};
		\draw [thick] plot[smooth, tension=.7] coordinates {(3.625,3.4) (3.25,2.75) (3.625,1.5) (2.875,0.625) (1.5,1) (0.75,0.5)};
		\node [draw,rotate=45,outer sep=0,inner sep=1,minimum size=12,fill=white] (v8) at (1.5,2.75) {$ b $};
		\node [draw,rotate=45,outer sep=0,inner sep=2,minimum size=12,fill=white] (v10) at (2.625,1.625) {$ b^* $};
		\node [draw,outer sep=0,inner sep=1,rotate=-45,minimum size=12,fill=white] (v12) at (3.75,3.75) {$i$};
		\node [draw,outer sep=0,inner sep=1,rotate=-45,minimum size=12,fill=white] (v11) at (0.5,0.5) {$j$};
		\node [draw,outer sep=0,inner sep=2,rotate=-45,minimum size=12,fill=white] (v5) at (-0.5,-0.5) {$k$};
		\node [draw,outer sep=0,inner sep=2,rotate=-45,minimum size=12,fill=white] (v13) at (-1.75,-1.75) {$l$};
		\draw[thick]  (v9) edge (v12);
		\draw [thick] (v5) edge (v11);
		\draw [thick] (v13) edge (v6);
		\draw [line width=0.18cm,white] (1.15,3.125) edge (v14);
		\draw [thick] (v8) edge (v14);
		\draw[very thick, red]  (v1) edge (v3);
		\draw[very thick, red]  (v2) edge (v4);
		\draw[very thick]  (v1) edge (v2);
		\draw[very thick]  (v3) edge (v4);
		\end{tikzpicture}
		\end{array} \\
		&\refequals[Lem.~\ref{LEMMA:DUAL_DECOMPOSE}] \frac{d(I)d(J)d(I')d(J')}{d(\C)^2}  \bigoplus\limits_{R} \sum\limits_{S} d(S)d(R)
		\begin{array}{c}
		\begin{tikzpicture}[scale=0.5,every node/.style={inner sep=0,outer sep=-1}]
		\node at (0,0) {};
		\node (v1) at (2.25,7.25) {};
		\node (v4) at (-0.25,-5.25) {};
		\node (v2) at (7.25,2.25) {};
		\node (v3) at (-5.25,-0.25) {};
		\node (v9) at (4.75,4.75) {};
		\node (v6) at (-2.75,-2.75) {};
		\node[outer sep=-0.5] (v7) at (4.625,-0.375) {};
		\node[outer sep=-0.5] (v14) at (-0.375,4.625) {};
		\node (v8) at (2.625,1.625) {};
		\node (v10) at (2.625,1.625) {};
		\node at (4,2.5) {$ J $};
		\node at (2.5,4) {$ I $};
		\node at (-1,1.5) {$ S $};
		\node at (1.5,-1) {$ S^\vee $};
		\node at (5.25,5.25) {$ X $};
		\node at (-3.25,-3.25) {$ Z $};
		\node at (-0.875,5.125) {$ R $};
		\node at (5.125,-0.875) {$ R $};
		\draw [thick] plot[smooth, tension=.7] coordinates {(3.5,3.75) (2.75,3.25) (1.75,3.5) (0.75,2.75) (1,1.5) (0.5,0.75)};
		\draw [thick] (v10) edge (v7);
		\draw [line width=0.18cm,white] (1.2,3.05) edge (v14);
		\draw [thick] (v8) edge (v14);
		\draw[thick]  plot[smooth, tension=.5] coordinates { (1.4128,1.7578) (2.125,0.375) (0,-1.625) (-0.625,-1.625)};
		\draw[line width=0.18cm,white]  plot[smooth, tension=.5] coordinates { (1.4128,1.7578) (0.375,2.125) (-1.625,0) (-1.625,-0.625)};
		\draw[thick]  plot[smooth, tension=.5] coordinates { (1.4128,1.7578) (0.375,2.125) (-1.625,0) (-1.625,-0.625)};
		\draw [line width=0.18cm,white] plot[smooth, tension=.7] coordinates { (3.25,2.75) (3.5,1.75) (2.75,0.75) (2.015,0.86) (1.5,1)};
		\draw [thick] plot[smooth, tension=.7] coordinates {(3.75,3.5) (3.25,2.75) (3.5,1.75) (2.75,0.75) (1.5,1) (0.75,0.5)};
		\draw [thick] plot[smooth, tension=.7] coordinates {(-1.625,-1.625) (-0.75,-1.5) (-0.625,-0.75)};
		\draw [thick] plot[smooth, tension=.7] coordinates {(-1.625,-1.625) (-1.625,-1.125) (-1.5522,-0.8822)};
		\draw [thick] plot[smooth, tension=.7] coordinates {(-0.75,-0.625) (-1.1,-0.6) (-1.3221,-0.6721)};
		\draw [thick] plot[smooth, tension=.7] coordinates {(-1.625,-0.625) (-0.875,-1.375)};
		\node [draw,outer sep=0,inner sep=1,rotate=-45,minimum size=12,fill=white] (v12) at (3.75,3.75) {$i$};
		\node [draw,outer sep=0,inner sep=1,rotate=-45,minimum size=12,fill=white] (v11) at (0.5,0.5) {$j$};
		\node [draw,outer sep=0,inner sep=2,rotate=-45,minimum size=12,fill=white] (v5) at (-0.5,-0.5) {$k$};
		\node [draw,outer sep=0,inner sep=2,rotate=-45,minimum size=12,fill=white] (v13) at (-1.75,-1.75) {$l$};
		\draw[thick]  (v9) edge (v12);
		\draw [thick] (v5) edge (v11);
		\draw [thick] (v13) edge (v6);
		\draw[very thick, red]  (v1) edge (v3);
		\draw[very thick, red]  (v2) edge (v4);
		\draw[very thick]  (v1) edge (v2);
		\draw[very thick]  (v3) edge (v4);
		\end{tikzpicture}
		\end{array}
	\end{align*}
	\begin{align*}
		\refequals[Prop.~\ref{PROP:TWISTED_S}] \delta_{I,I'} \delta_{J,J'} \tr(k \circ jr) \frac{d(I)d(J)}{d(\C)}\bigoplus\limits_{R} d(R)		\begin{array}{c}
		\begin{tikzpicture}[scale=0.5,every node/.style={inner sep=0,outer sep=-1}]
		\node (v1) at (0,5) {};
		\node (v4) at (0,-5) {};
		\node (v2) at (5,0) {};
		\node (v3) at (-5,0) {};
		\node (v9) at (2.5,2.5) {};
		\node (v6) at (-2.5,-2.5) {};
		\node (v8) at (-2.5,2.5) {};
		\node (v11) at (2.5,-2.5) {};
		\node [every node] (v10) at (0.8,-0.8) {};
		\node [every node] (v12) at (1.2,-1.2) {};
		\node [every node] (v13) at (-1.5,-1.25) {};
		\node [every node] (v14) at (-1.2,0.8) {};
		\node [every node] (v15) at (-0.8,1.2) {};
		\node (v20) at (1.25,1.5) {};
		\node (v16) at (-1.25,-1.5) {};
		\node [every node] (v18) at (1,-1) {};
		\node (v17) at (1.5,1.25) {};
		\node [every node] at (3,3) {$ X $};
		\node [every node] at (-2.0068,-0.0314) {$ I $};
		\node [every node] at (1.9687,-0.0314) {$ J $};
		\node [every node] at (-3,-3) {$ L $};
		\node [every node] at (-3,3) {$ R $};
		\node [every node] at (3,-3) {$ R $};
		\draw [thick] (v8) edge (v10);
		\draw [thick] (v11) edge (v12);
		\draw[thick]  plot[smooth, tension=.7] coordinates {(v20) (0.6,1) (-0.2,1.4) (v15)};
		\draw[thick]  plot[smooth, tension=.7] coordinates {(v13) (-1,-0.6) (-1.4,0.2) (v14)};
		\draw[thick]  plot[smooth, tension=.7] coordinates {(v16) (-0.6,-1) (0.2,-1.4) (v18) (1.4,-0.2) (1,0.6) (v17)};
		\node [draw,rotate=-45,outer sep=0,inner sep=1,minimum size=12,fill=white] (v5) at (1.5,1.5) {$i$};
		\node [draw,rotate=-45,outer sep=0,inner sep=2,minimum size=12,fill=white] (v7) at (-1.5,-1.5) {$l$};
		\draw[thick]  (v9) edge (v5);
		\draw [thick] (v6) edge (v7);
		\draw[very thick, red]  (v1) edge (v3);
		\draw[very thick, red]  (v2) edge (v4);
		\draw[very thick]  (v1) edge (v2);
		\draw[very thick]  (v3) edge (v4);
		\end{tikzpicture}
		\end{array}.
		\end{align*}
	\end{REM}
	\bibliographystyle{alpha}
	\bibliography{Paper_Version}

\begin{thebibliography}{EGNO15}

\bibitem[BK01]{GraphCal}
Bojko Bakalov and Alexander Kirillov, Jr.
\newblock {\em Lectures on tensor categories and modular functors}, volume~21
  of {\em University Lecture Series}.
\newblock American Mathematical Society, Providence, RI, 2001.

\bibitem[EGNO15]{EtingofBook}
Pavel Etingof, Shlomo Gelaki, Dmitri Nikshych, and Victor Ostrik.
\newblock {\em Tensor categories}, volume 205 of {\em Mathematical Surveys and
  Monographs}.
\newblock American Mathematical Society, Providence, RI, 2015.

\bibitem[EK98]{Evans1998}
David~E. Evans and Yasuyuki Kawahigashi.
\newblock Orbifold subfactors from {Hecke Algebras II} -- quantum doubles and
  braiding.
\newblock {\em Comm. Math. Phys.}, 196(2):331--361, Aug 1998.

\bibitem[Kon08]{MR2430629}
Liang Kong.
\newblock Cardy condition for open-closed field algebras.
\newblock {\em Comm. Math. Phys.}, 283(1):25--92, 2008.

\bibitem[M\"{u}3]{MR1990929}
Michael M\"{u}ger.
\newblock On the structure of modular categories.
\newblock {\em Proc. London Math. Soc. (3)}, 87(2):291--308, 2003.

\bibitem[ML98]{MR1712872}
Saunders Mac~Lane.
\newblock {\em Categories for the working mathematician}, volume~5 of {\em
  Graduate Texts in Mathematics}.
\newblock Springer-Verlag, New York, second edition, 1998.

\bibitem[Ocn94]{Ocneanu1993}
Adrian Ocneanu.
\newblock Chirality for operator algebras.
\newblock In {\em Subfactors ({K}yuzeso, 1993)}, pages 39--63. World Sci.
  Publ., River Edge, NJ, 1994.

\bibitem[PSV18]{Vaes}
Sorin Popa, Dimitri Shlyakhtenko, and Stefaan Vaes.
\newblock Cohomology and {$L^2$}-{B}etti numbers for subfactors and
  quasi-regular inclusions.
\newblock {\em Int. Math. Res. Not. IMRN}, (8):2241--2331, 2018.

\end{thebibliography}
\end{document}